# An Upwind Generalized Finite Difference Method (GFDM) for Meshless Analysis of Heat and Mass Transfer in Porous Media

Xiang Rao*

[1]Cooperative Innovation Center of Unconventional Oil and Gas (Ministry of Education & Hubei Province), Yangtze University, Wuhan, 430100, China

[2]Key Laboratory of Drilling and Production Engineering for Oil and Gas, Hubei Province, Wuhan 430100, Hubei, China

[3]School of Petroleum Engineering, Yangtze University, Wuhan 430100, China

*Corresponding author: raoxiang0103@163.com, raoxiang@yangtzeu.edu.cn

**Abstract**

In this paper, an upwind GFDM is developed for the coupled heat and mass transfer problems in porous media. GFDM is a meshless method that can obtain the difference schemes of spatial derivatives by using Taylor expansion in local node influence domains and the weighted least squares method. The first-order single-point upstream scheme in the FDM/FVM-based reservoir simulator is introduced to GFDM to form the upwind GFDM, based on which, a sequential coupled discrete scheme of the pressure diffusion equation and the heat convection-conduction equation is solved to obtain pressure and temperature profiles. This paper demonstrates that this method can be used to obtain the meshless solution of the convection-diffusion equation with a stable upwind effect. For porous flow problems, the upwind GFDM is more practical and stable than the method of manually adjusting the influence domain based on the prior information of the flow field to achieve the upwind effect. Two types of calculation errors are analyzed, and three numerical examples are implemented to illustrate the good calculation accuracy and convergence of the upwind GFDM for heat and mass transfer problems in porous media, and indicate the increase of the radius of the node influence domain will increase the calculation error of temperature profiles. Overall, the upwind GFDM discretizes the computational domain using only a point cloud that is generated with much less topological constraints than the generated mesh, but achieves good computational performance as the mesh-based approaches, and therefore has great potential to be developed as a general-purpose numerical simulator for various porous flow problems in domains with complex geometry.

## 1. Introduction

The study on heat and mass transfer in porous media widely exists in the development and utilization of environment-friendly geothermal resources, thermal recovery of oil and gas resources, thermal performance of insulation materials, etc. Underground formations are typical porous media, among which the study on coupled heat and mass transfer focuses on the coupling calculation of fluid seepage and heat conduction-convection in porous formation. This coupling effect is mainly reflected in the aspects [1], including the influence of temperature change on fluid viscosity, the influence of temperature on formation porosity, and the influence of fluid flow velocity on the strength of heat convection, etc. At present, the numerical simulation methods of coupled mass and heat transfer mainly include finite difference method (FDM) [2], finite element method (FEM) [3], and finite volume method (FVM) [4, 5]. However, these methods are limited by the requirements of geometric regularity of the computational domain and high-quality mesh generation.

The generalized finite difference method (GFDM) is a domain-type meshless method with twenty years. Based on Taylor series expansion of unknown function and weighted least square approximation in a subdomain, the spatial derivatives of unknown function in the governing equation are expressed as the difference scheme of the values of unknown function at nodes in the subdomain, which overcomes the grid dependence of traditional FDM [6, 7]. Up to now, GFDM has been widely used to solve various scientific and engineering problems, including coupled thermoelasticity problem [8-10], third and fourth-order partial differential equations [11], shallow water equations [12], transient heat conduction analysis [13], seismic wave propagation problem [14], stress analysis [15], unsteady Burgers' equations [16-17], water wave interactions [18], inverse heat source problems [19], nonlinear convection-diffusion equations [20], time-fractional diffusion equation [21], various flow problems [22, 23, 24]. Gavete et al. [25] reviewed the advantages and disadvantages of GFDM and its applications, analyzed the influence of various factors on the numerical performances of GFDM, and found that the weight function might have little influence on the numerical results. GFDM just uses point clouds for the discretization of the computational domain to realize the effective numerical solution of partial equations, which saves the possible time-consuming and laborious meshing and numerical integration in FEM, FVM, and boundary element method (BEM) [26, 27] for the calculation domains with complex geometry.

This paper focuses on applying GFDM to the modeling of coupled heat and mass transfer problems, which involves not only the diffusion equation about pressure, but also the convection-diffusion equation about temperature. For the convection-diffusion equation, it is often necessary to add upwind weight treatment to the discrete scheme of the convection term, otherwise, the calculation solution is prone to the situation of inaccurate oscillation. For example, the upstream FEM [28], the upstream FDM [29], and the upstream FVM [30] have been widely used. In meshless methods, currently, modifying the influence domain is generally adopted to realize the upwind effect, including the upwind influence domain [31] of moving the central node position in the upstream direction and the partial influence domain [32, 12, 33] of including the upstream nodes more in the central-node influence domain. For GFDM, Cheng and Liu [32] roughly discussed the upwind effect by constructing a six-point scheme containing more upstream nodes in the influence domain (i.e., the partial influence domain) in GFDM. Li and Fan [12] adopted the partial influence domain to handle the convection-dominated hyperbolic shallow water equations, then used the flux limiter technique to avoid the non-physical wiggles of solutions near discontinuities [33]. However, for porous flow problems, including the coupled calculation of mass and heat transfer in porous media studied in this paper, because the underground velocity field is generally unknown, and the velocity field of reservoir flow may be very complex due to the influence of various geological conditions or source-sink terms caused by drilling wells, it is difficult to obtain a stable upwind effect by modifying the influence domain to ensure good calculation performance.

Therefore, this paper aims to study the computational performance of the upwind GFDM to heat and mass transfer problem, including the convection-dominated convection-diffusion equation, so as to provide an important reference for constructing a general-purpose numerical simulator for multiphysics coupling porous flow problems based on the upwind GFDM.

The paper is structured as follows. In Section 2, the upwind GFDM based modeling of the single-phase heat and mass transfer problem is given, including the basic physical model in Section 2.1, a brief review of GFDM in Section 2.2, the upwind GFDM discrete schemes of heat and mass transfer equations in Section 2.3, the treatment of boundary conditions in Section 2.4, and the application of the upwind GFDM to the meshless solution of the convection-diffusion equation and analysis of the dissipation error in Section 2.5. Section 3 gives three numerical examples and a rough error analysis to illustrate the computational performances of the upwind GFDM. The conclusion and future work come in Section 4.

Nomenclature

| Symbols | Physical meanings |
|---|---|
| $\alpha$ | the unit conversion factor, equal to 0.0864 |
| $\beta$ | the unit conversion factor, equal to 86400 |
| $\bar{v}$ | the seepage velocity, m/day |
| $k$ | permeability, mD |
| $\mu$ | viscosity, cp |
| $p$ | pressure, MPa |
| $T$ | temperature, °C |
| $q$ | source or sink term in mass transfer, 1/day |
| $\phi$ | formation porosity which is a function of pressure and temperature $\phi(p,T)$, fraction |
| $t$ | time, day |
| $\lambda_c$ | integrated heat conduction coefficient, a function of pressure and temperature $\lambda_c(p,T)$, J/s/m/°C |
| $\lambda_l$ | heat conduction coefficient of liquid, J/s/m/°C |
| $\lambda_r$ | heat conduction coefficient of rock, J/s/m/°C |
| $\rho_l$ | liquid density, kg/m³ |
| $\rho_r$ | rock density, kg/m$^3$ |
| $C_l$ | liquid heat capacity, J/kg/°C |
| $C_r$ | rock heat capacity, J/kg/°C |
| $q_h$ | the energy source or sink term which is the total heat energy carried by the mass source or sink term $q$, J/m$^3$/day |
| $C_t$ | compressibility coefficient, 1/MPa |
| $C_{Temp}$ | the thermal expansion coefficient, 1/°C |
| $p_0$ | the initial formation pressure, MP |
| $T_0$ | the initial formation temperature, °C |
| $\phi_0$ | the porosity when $p=p_0$ and $T=T_0$ |
| $\alpha_T$ | the fluid viscosity-temperature coefficient, which measures the physical law that the fluid viscosity decreases with the increase of the temperature |

## 2. Methodology
### 2.1 Governing equations

This paper focuses on the study of single-phase heat and mass transfer in porous media, including diffusion equation about pressure, convection-diffusion equation about temperature, and auxiliary equations of physical quantities affected by pressure and temperature, which are:

(1) Mass conservative equation (Assuming that the fluid is incompressible)

$$-\nabla \cdot \vec{v} + q = \frac{\partial \phi(p,T)}{\partial t} \tag{1}$$

(2) Energy conservative equation

$$\beta\nabla\left(\lambda_c(p,T)\nabla T\right) - \nabla\cdot\left(\rho_l C_l T \vec{v}\right) + q_h = \frac{\partial}{\partial t}\left(\left[1-\phi(p,T)\right]\rho_r C_r T + \phi(p,T)\rho_l C_l T\right) \tag{2}$$

(3) Auxiliary equations

In porous flow, the seepage velocity $\vec{v}$ in Eq. (1) and Eq. (2) satisfies Darcy's law:

$$\vec{v} = -\alpha \frac{k}{\mu}\nabla p \tag{3}$$

Thus, Eq. (1) is rewritten as a form of an approximate pressure diffusion equation:

$$\alpha\nabla\left(\frac{k}{\mu}\nabla p\right) + q = \frac{\partial \phi(p,T)}{\partial t} \tag{4}$$

Eq. (2) can be rewritten as:

$$\beta\nabla\left(\lambda_c(p,T)\nabla T\right) + \alpha\nabla\left(\rho_l C_l T \frac{k}{\mu}\nabla p\right) + q_h = \frac{\partial}{\partial t}\left(\left[1-\phi(p,T)\right]\rho_r C_r T + \phi(p,T)\rho_l C_l T\right) \tag{5}$$

Due to the elastic and thermoelastic properties of reservoir porous media, porosity is affected by both fluid pressure and temperature. The porosity $\phi(p,T)$, integrated heat conduction coefficient $\lambda_c(p,T)$, and liquid viscosity $\mu(T)$ are calculated as

$$\begin{aligned}\phi(p,T) &= \left[\phi_0 + C_t(p-p_0)\right]\left[1 + \frac{1-\phi_0}{\phi_0}C_{Temp}(T-T_0)\right], \\ \lambda_c(p,T) &= \phi(p,T)\lambda_l + (1-\phi(p,T))\lambda_r, \\ \mu(T) &= \mu(T_0)e^{-\alpha_T(T-T_0)}\end{aligned} \tag{6}$$

Eq. (4) and Eq. (5) are differential equations about temperature and pressure, in which the coefficients are jointly affected by pressure and temperature to form a nonlinear equation system.

### 2.2 A brief review of GFDM

GFDM is a relatively new meshless method based on local Taylor expansion and weighted least squares approximation. In this method, the spatial derivatives are approximated as a difference scheme of the nodal function values within each local node influence domain.

Suppose there are $n$ other nodes in the influence domain of the node $M_0 = (x_0, y_0)$, which are denoted as $\{M_1, M_2, M_3, \cdots, M_N\}$ where $M_i = (x_i, y_i)$. The Taylor expansion of the unknown function $u(x,y)$ at $M_0$ can be used to approximate $\{u(M_i), i=1,\cdots n\}$ as:

$$u(M_i) = u(M_0) + \Delta x_i \left.\frac{\partial u}{\partial x}\right|_{M_0} + \Delta y_i \left.\frac{\partial u}{\partial y}\right|_{M_0} + \frac{1}{2}\left((\Delta x_i)^2 \left.\frac{\partial^2 u}{\partial x^2}\right|_{M_0} + 2\Delta x_i \Delta y_i \left.\frac{\partial^2 u}{\partial x \partial y}\right|_{M_0} + (\Delta y_i)^2 \left.\frac{\partial^2 u}{\partial y^2}\right|_{M_0}\right) + O(r^3) \quad (7)$$

where $\Delta x_i = x_0 - x_i$, $\Delta y_i = y_0 - y_i$.

Denote $u_0 = u(M_0)$, $u_{x0} = \left.\frac{\partial u}{\partial x}\right|_{M_0}$, $u_{y0} = \left.\frac{\partial u}{\partial y}\right|_{M_0}$, $u_{xx0} = \left.\frac{\partial^2 u}{\partial x^2}\right|_{M_0}$, $u_{xy} = \left.\frac{\partial^2 u}{\partial x \partial y}\right|_{M_0}$, $u_{yy} = \left.\frac{\partial^2 u}{\partial y^2}\right|_{M_0}$.

Define weighted error function $B(\mathbf{D}_u)$:

$$B(\mathbf{D}_u) = \sum_{j=1}^{n}\left[\left(u_0 - u_j + \Delta x_j u_{x0} + \Delta y_j u_{y0} + \frac{1}{2}(\Delta x_j)^2 u_{xx0} + \frac{1}{2}(\Delta y_j)^2 u_{yy0} + \Delta x_j \Delta y_j u_{xy0}\right)\omega_j\right]^2 \quad (8)$$

where $\mathbf{D}_u = (u_{x0}, u_{y0}, u_{xx0}, u_{yy0}, u_{xy0})^T$, $\omega_j = \omega(\Delta x_j, \Delta y_j)$ is the value of the weight function $\omega(x,y)$ at $M_j$, Benito et al. [6] and Gavete et al. [25] demonstrated that different types of weight function have little influence on the calculation results, while the quartic spline function is generally selected as the weight function in Eq. (9).

$$\omega_j = \begin{cases} 1 - 6\left(\frac{r_j}{r_m}\right)^2 + 8\left(\frac{r_j}{r_m}\right)^3 - 3\left(\frac{r_j}{r_m}\right)^4 & r_j \leq r_m \\ 0 & r_j > r_m \end{cases} \quad (9)$$

where $r_j$ is the Euclidean distance from the $M_j$ to $M_0$, and $r_m$ is the radius of the influence domain of $M_0$.

The weighted error function $B(\mathbf{D}_u)$ is minimized, at this time, the partial derivatives of $B(\mathbf{D}_u)$ to each component of $B(\mathbf{D}_u)$ are required equal to zero, they are

$$\frac{\partial B(\mathbf{D}_u)}{\partial u_{x0}} = 0, \frac{\partial B(\mathbf{D}_u)}{\partial u_{y0}} = 0, \frac{\partial B(\mathbf{D}_u)}{\partial u_{xx0}} = 0, \frac{\partial B(\mathbf{D}_u)}{\partial u_{yy0}} = 0, \frac{\partial B(\mathbf{D}_u)}{\partial u_{xy0}} = 0 \quad (10)$$

The above equations in Eq. (9) are sorted into linear equations as follows:

$$\mathbf{A}\mathbf{D}_u = \mathbf{b} \quad (11)$$

where $\mathbf{A} = \mathbf{L}^T \boldsymbol{\omega} \mathbf{L}$, $\mathbf{b} = \mathbf{L}^T \boldsymbol{\omega} \mathbf{U}$, $\mathbf{L} = (\mathbf{L}_1^T, \mathbf{L}_2^T, \cdots, \mathbf{L}_n^T)^T$, $\mathbf{L}_i = \left(\Delta x_i, \Delta y_i, \frac{\Delta x_i^2}{2}, \frac{\Delta y_i^2}{2}, \Delta x_i \Delta y_i\right)$, $\boldsymbol{\omega} = diag(\omega_1^2, \omega_2^2, \cdots, \omega_n^2)$, $\mathbf{U} = (u_1 - u_0, u_2 - u_0, \cdots, u_n - u_0)^T$.

Then, $\mathbf{D}_u$ can be solved as:

$$\mathbf{D}_u = (u_{x0}, u_{y0}, u_{xx0}, u_{yy0}, u_{xy0})^T = \mathbf{A}^{-1}\mathbf{b} = \mathbf{A}^{-1}\mathbf{L}^T \boldsymbol{\omega} \mathbf{U} = \mathbf{M}\mathbf{U} \quad (12)$$

where $\mathbf{M} = \mathbf{A}^{-1}\mathbf{L}^T \boldsymbol{\omega}$。

For the convenience of notation, the elements of the matrix $\mathbf{M}$ are denoted as $m_{ij}$, and the generalized difference approximation schemes of the spatial derivatives at $M_0$ are obtained as:

$$\frac{\partial u}{\partial x} = \sum_{j=1}^{n} m_{1j}(u_j - u_0), \quad \frac{\partial u}{\partial y} = \sum_{j=1}^{n} m_{2j}(u_j - u_0), \quad \frac{\partial^2 u}{\partial x^2} = \sum_{j=1}^{n} m_{3j}(u_j - u_0), \quad \frac{\partial^2 u}{\partial y^2} = \sum_{j=1}^{n} m_{4j}(u_j - u_0), \quad (13)$$

$$\frac{\partial^2 u}{\partial x \partial y} = \sum_{j=1}^{n} m_{5j}\left(u_j - u_0\right)$$

As seen in Eq. (13), GFDM is flexible to obtain the difference expressions of the first-order and second-order spatial derivatives at the considered node only according to the coordinates of the nodes within the influence domain of the considered node. In fact, there can be no concept of the node influence domain, because it is only necessary to determine which nodes participate in the construction of the generalized finite difference expressions of spatial derivatives at the considered node, and the introduction of the node influence domain is just to determine the selection of these nodes more conveniently.

Therefore, the numerical discretization of partial differential equations can be realized when only using a point cloud to discretize the computational domain. This is the most significant advantage of meshless GFDM compared with mesh-based FEM and FVM.

Milewski [35] and Rao et al. [34] point out that the point cloud discretization in the computational domain has much less topological information than the mesh discretization in the computational domain, for example, on the basis of the point cloud, the mesh also needs to determine which two points are connected to an edge, which points form a mesh and the order of the vertices of the mesh, and so on. When generating a mesh, the lengths of the edges of a mesh do not vary so much that the vertex angles of the mesh do not vary so much to ensure the quality of the mesh generation. Therefore, when discretizing a computational domain, the topological constraints on the mesh generation are much greater than those on the generation of the point cloud, which makes the generation of the point cloud theoretically much less difficult than the mesh division.

Rao et al. [34] takes a circular domain as an example, and points out that when the point cloud is used to discretize the computing domain, the method of equally dividing the radius and argument can be quickly used to form a concentric point cloud, and the workload of meshing the computing domain is obviously larger than that of generating the concentric point cloud. Even if the mesh is constructed on the basis of this point cloud, it is difficult to give a criterion for which nodes in the point cloud form the mesh to determine the high-quality mesh generation. This shows that the discretization of point cloud to computing domain can be more arbitrary and simpler than that of grid to computing domain, which is precisely because the generation of the point cloud is limited by little topology.

In terms of point cloud generation methods for computational domains, Milewski [35] showed that Liszka-type node generators, which were proposed by Liszka [31], can be used to generate point clouds that can be adequately adapted to irregular computational domains. Michel et al. [23] applied a meshfree advancing front technique [32] to generate the initial point cloud. Löhner and Onate [37] developed an algorithm to construct boundary-conforming, isotropic clouds of points with variable density in space, which is more efficient than mesh-generation methods to adaptively discretize the computational domain. Rao et al. [38] gave an algorithm for point cloud generation for computational domains with complex geometrical entities in lower one dimension, such as complex fracture networks in fractured reservoirs, showing that point cloud discretization of computational domains can effectively solve the gridding challenge of matching grids for fracture intersections and narrow computational domains between fractures, avoiding the problem of generating very fine grids in these locally narrow areas due to the length of each side of the grid being required to differ significantly. Moreover, the algorithm for generating grid nodes in the grid division method can also be used directly to generate point clouds of the computational domain, for example, Cartesian collocation points or mesh vertices of the triangulation can be used as point clouds of the computational domain. Overall, the generation of point clouds is less difficult than mesh division, and the generation methods are more diverse and easier to carry out an adaptive analysis. The above analysis demonstrates the advantages of meshless methods such as GFDM over mesh-based numerical methods in terms of computational domain discretization, and this was an important original motivation for the rapid development of meshless methods in the 1990s.

**2.3 Upwind GFDM based discrete schemes**

After the calculation domain has been discretized by a point cloud, the nodes in the point cloud are denote as node $i$, $i = 1, 2, 3, \cdots, n_t$, in which $n_t$ is the total number of nodes. Define the set composed of the sequence numbers of the nodes in the influence domain of node $i$ except the node $i$ itself (i.e. the nodes participating in the construction of the generalized difference operator of node $i$ except the node $i$ itself) as the index set of node $i$, which is denoted as $\Lambda_i$. if the influence domain of node $i$ contains $n_i$ nodes except node $i$, then $\Lambda_i$ has $n_i$ elements. The nodes in $\Lambda_i$ and the node $i$ itself together form the local point cloud of node $i$.

For node $i$, according to Eq. (13) in GFDM, it is obtained that:

$$\frac{\partial^2 p}{\partial x^2} = \sum_{j \in \Lambda_i} m_{3j}^i (p_j - p_i), \quad \frac{\partial^2 p}{\partial y^2} = \sum_{j \in \Lambda_i} m_{4j}^i (p_j - p_i) \tag{14}$$

where the superscript $i$ of $m_{3j}^i$ indicates that node $i$ is the considered node.

Then, the pressure diffusion term is approximated as follows:

$$\nabla \cdot (\nabla p) = \frac{\partial^2 p}{\partial x^2} + \frac{\partial^2 p}{\partial y^2} \approx \sum_{j \in \Lambda_i} \left[ \left( m_{3j}^i + m_{4j}^i \right) (p_j - p_i) \right] \tag{15}$$

For the actual underground formation, permeability $k$ is often difficult to express as an explicit function about coordinates, but only knows the permeability values at some nodes, therefore, it is difficult for us to take $k$ as a function and extract it from the diffusion term. Therefore, this paper uses the harmonic average scheme of the nodal permeability values to calculate the permeability between node $i$ and node $j$, the arithmetic average scheme of the nodal viscosity values is used to characterize the fluid viscosity between node $i$ and node $j$, which are commonly used in FDM/FVM-based reservoir simulator (i.e., the numerical simulator about porous flow problems) [38, 39, 40, 41]. The treatments of the heterogeneity of physical parameters are beneficial to the easier application of GFDM in practical porous flow problems, because it is generally difficult to obtain the function expression of physical parameters with good smoothness in practical problems, especially the related physical parameters of underground reservoirs. The numerical examples in Section 3.2 will prove that such treatment can achieve sufficient calculation accuracy. Therefore, it is obtained that:

$$\nabla \cdot \left( \frac{k}{\mu(T)} \nabla p \right) = \sum_{j \in \Lambda_i} \left[ \frac{k_{ij}}{\mu_{ij}} \left( m_{3j}^i + m_{4j}^i \right) (p_j - p_i) \right] \tag{16}$$

where $k_{ij} = \dfrac{2 k_i k_j}{k_i + k_j}$, $\mu_{ij} = \dfrac{\mu(T_i) + \mu(T_j)}{2}$.

For the heat convection-diffusion equation, the heat conduction term in Eq. (2) can adopt a discrete scheme similar to Eq. (15), and obtain:

$$\nabla (\lambda_c (p,T) \nabla T) = \sum_{j \in \Lambda_i} \left[ \lambda_{c,ij} \left( m_{3j}^i + m_{4j}^i \right) (T_j - T_i) \right] \tag{17}$$

where $\lambda_{c,ij} = \dfrac{2 \lambda_{c,i} \lambda_{c,j}}{\lambda_{c,i} + \lambda_{c,j}}$, $\lambda_{c,i} = \lambda_c (p_i, T_i)$, $\lambda_{c,j} = \lambda_c (p_j, T_j)$.

For the heat convection term in Eq. (2), because the convection term has asymmetry, its discretization needs to take the upwind scheme. In a meshless method, the upstream effect is generally constructed by modifying the influence domain, however, this method is not conducive to taking a stable upwind effect with a complex flow field and the construction of a general framework, therefore, this paper aims to constructs the upstream scheme in GFDM for porous flow problems without modifying the influence domain.

In the porous flow problem, the velocity in porous media satisfies Darcy's law in Eq. (3), then in the porous flow problem, the heat convection term $-\rho_l C_l \nabla \cdot (\vec{v} T)$ (It should be $-\nabla \cdot (\rho_l C_l \vec{v} T)$, but considering that $\rho_l$ and $C_l$ can be regarded as constant, so for the convenience of discussion and analysis, we move them to the outside of the Hamiltonian) in Eq. (2) can be rewritten as Eq. (18) in Eq. (5).

$$-\rho_l C_l \nabla \cdot (T \vec{v}) = \alpha \rho_l C_l \nabla \left( T \frac{k}{\mu} \nabla p \right) \tag{18}$$

Thus obtain a second-order derivative form similar to the diffusion term. The FDM/FVM-based reservoir simulator generally uses the first-order single point upwind (SPU) scheme to discretize the convection term. Taking the FDM as an example, if the difference scheme of pressure diffusion term $\nabla \left( \dfrac{k}{\mu} \nabla p \right)$ is:

$$\nabla\left(\frac{k}{\mu}\nabla p\right)=\sum_{j}D_{ij}\left(p_{j}-p_{i}\right) \tag{19}$$

where $D_{ij}$ is the coefficient in the FDM-based difference expression.

Then $D_{ij}\left(p_{j}-p_{i}\right)$ nearly denote (not rigorously) the seepage velocity (information) of the fluid between grid $i$ and grid $j$, and the difference scheme of convection term $\alpha\rho_{l}C_{l}\nabla\left(T\frac{k}{\mu}\nabla p\right)$ is discretized as:

$$\alpha\rho_{l}C_{l}\nabla\left(T\frac{k}{\mu}\nabla p\right)=\alpha\rho_{l}C_{l}\sum_{j}T_{ij}D_{ij}\left(p_{j}-p_{i}\right) \tag{20}$$

where, if node $i$ is the upstream of node $j$, that is, if $p_i > p_j$, the SPU scheme obtains $T_{ij}=T_i$, otherwise, $T_{ij}=T_j$, then,

$$\alpha\rho_{l}C_{l}T_{ij}D_{ij}\left(p_{j}-p_{i}\right)=\alpha\rho_{l}C_{l}T_{i}D_{ij}\left(p_{j}-p_{i}\right) \tag{21}$$

roughly represents the heat loss in grid $i$ caused by the flow of upstream grid $i$ to downstream grid $j$. It can be seen that for the convection term in the porous flow problem controlled by the pressure gradient, the difference expression of the pressure diffusion term in Eq. (19) can be obtained first, which contains the seepage velocity information between the central grid and the adjacent grid, and then the SPU scheme is adopted for the physical quantities related to convection transfer (such as the temperature $T$ in the thermal convection term) to realize the discretization of the convection term.

The SPU scheme in FDM/FVM-based reservoir simulator provides a great inspiration for the discretization of the convection term when the meshless method is applied to the porous flow problems, because it is simple to obtain the difference scheme of the pressure diffusion term by using the meshless method, and it seems that the SPG scheme in Eq. (21) can also be directly applied to the meshless difference scheme by replacing the pressure of grid $i$ and grid $j$ with the pressure of node $i$ and node $j$, The work done in this paper is to verify whether the SPG scheme in GFDM can achieve good calculation performance for mass and heat transfer in porous media, so as to form an upwind GFDM used for porous flow problems which can achieve a stable upwind effect.

Therefore, the SPU scheme is directly introduced to GFDM to form an upwind GFDM. Inspired by this idea, because the difference scheme of pressure diffusion term $\nabla\left(\frac{k}{\mu}\nabla p\right)$ in GFDM is Eq. (16), the discrete scheme in Eq. (22) of the heat convection term with the upwind scheme of temperature in Eq. (23) is adopted. In Section 2.5, we will demonstrate that when dealing with the one-dimensional constant-coefficient heat convection-diffusion equation, the discretization of the heat convection term by the upwind GFDM can be reduced to the discretization of the convection term by FDM with first-order upwind scheme. The computational performance of the upwind GFDM will be verified and analyzed by numerical examples in Section 3.

$$\nabla\left(\rho C_{l}T\frac{k}{\mu}\nabla p\right)=\sum_{j\in\Lambda_{i}}\left[\rho C_{l}T_{ij}\frac{k_{ij}}{\mu_{ij}}\left(m_{3j}^{i}+m_{4j}^{i}\right)\left(p_{j}-p_{i}\right)\right] \tag{22}$$

$$T_{ij}=\begin{cases}T_{j} & \text{if } p_{j}\geq p_{i}\\ T_{i} & \text{if } p_{j}<p_{i}\end{cases} \tag{23}$$

For the convection term expressed by the first-order spatial derivative in the hyperbolic shallow water equation, if it can be rewritten into the form expressed by the second-order spatial derivative similar to the diffusion term in porous-media seepage mechanics, it may also be solved by the upwind GFDM given in this paper. If not, the developed upwind GFDM may not be applicable to the high-performance solution of the shallow water equation, however, this method should be suitable for solving various flow problems in porous media, which is also the original motivation of this paper, and its good performance in hyperbolic two-phase flow in porous media has been confirmed [34].

Then, a sequential coupled scheme is adopted, that is, based on the temperature values of *n* time steps, implicitly calculate the pressure values of *n*+1 time step, and then calculate the temperature values of *n*+1 time step.

Therefore, the discrete scheme of the right side of the pressure diffusion equation is as follows:

$$\frac{\partial \phi(p,T)}{\partial t} = \frac{\partial \phi(p,T^n)}{\partial t} = \left[1 + \frac{1-\phi_0}{\phi_0} C_{Temp}(T^n - T_0)\right] C_t \frac{\partial p}{\partial t} = \left[1 + \frac{1-\phi_0}{\phi_0} C_{Temp}(T^n - T_0)\right] C_t \frac{p^{n+1} - p^n}{\Delta t} \quad (24)$$

where $\Delta t$ is the time interval between the *n* time step and the *n*+1 time step.

Then the discrete scheme of Eq. (4) is obtained as follows:

$$\alpha \sum_{j \in \Lambda_i} \frac{k_{ij}}{\mu_{ij}}\left[\left(m_{3j}^i + m_{4j}^i\right)\left(p_j^{n+1} - p_i^{n+1}\right)\right] + q_i^{n+1} = \left[1 + \frac{1-\phi_0}{\phi_0} C_{Temp}(T_i^n - T_0)\right] C_t \frac{p_i^{n+1} - p_i^n}{\Delta t} \quad (25)$$

When the source or sink term is zero, the linear equation is sorted as follows:

$$\alpha \sum_{j \in \Lambda_i} \frac{k_{ij}}{\mu_{ij}^n}\left(m_{3j}^i + m_{4j}^i\right) p_j^{n+1} - \left(\alpha \sum_{j \in \Lambda_i} \frac{k_{ij}}{\mu_{ij}^n}\left(m_{3j}^i + m_{4j}^i\right) + \frac{1}{\Delta t}\left[1 + \frac{1-\phi_0}{\phi_0} C_{Temp}(T_i^n - T_0)\right]\right) C_t p_i^{n+1}$$

$$+ q_i^{n+1} = -\frac{1}{\Delta t}\left[1 + \frac{1-\phi_0}{\phi_0} C_{Temp}(T_i^n - T_0)\right] C_t p_i^n \quad (26)$$

where $\mu_{ij}^n = \dfrac{\mu(T_i^n) + \mu(T_j^n)}{2}$.

By synthesizing the discrete pressure diffusion equation at each node, combined with the boundary condition (the treatments of boundary condition will be illustrated in Section 2.4), global linear equations can be obtained to solve nodal pressure values at *n*+1 time step. Then, the discrete equations of Eq. (5) about temperature distribution are obtained:

$$\beta \sum_{j \in \Lambda_i}\left[\lambda_{c,ij}(p^{n+1}, T^n)(m_{3j}^i + m_{4j}^i)(T_j^{n+1} - T_i^{n+1})\right] + \alpha \sum_{j \in \Lambda_i}\left[\rho_l C_l T_{ij}^{n+1} \frac{k_{ij}}{\mu_{ij}^n}(m_{3j}^i + m_{4j}^i)(p_j^{n+1} - p_i^{n+1})\right]$$

$$+ q_H^{n+1} = \frac{\left[(1-\phi_i^{n+1})\rho_r C_r + \phi_i^{n+1}\rho_l C_l\right] T_i^{n+1} - \left[(1-\phi_i^n)\rho_r C_r + \phi_i^n \rho_l C_l\right] T_i^n}{\Delta t} \quad (27)$$

where $\phi_i^{n+1} = \phi(p_i^{n+1}, T_i^n)$, $\phi_i^n = \phi(p_i^n, T_i^n)$, $\mu_{ij}^n = \dfrac{\mu(T_i^n) + \mu(T_j^n)}{2}$, $T_{ij} = \begin{cases} T_j & \text{if } p_j \geq p_i \\ T_i & \text{if } p_j < p_i \end{cases}$.

Due to the sequential coupling scheme, the time step will be relatively small. Therefore, in the actual calculation, $T_{ij}^{n+1}$ in the discretization of the convection term in Eq. (27) can also be taken as $T_{ij}^n$.

## 2.4 Treatment of boundary conditions

As can be seen in Section 2.2, GFDM employs node coordinate information in the node influence domain in conjunction with Taylor expansion to obtain the generalized difference operator that minimizes the weighted truncation errors of Taylor expansion, i.e., the weighted least square method commonly used to obtain the spatial derivatives of unknown functions in various meshless methods. As a result, the quality of the node distribution participating in the construction of the generalized finite difference expression of spatial derivatives at the considered node significantly affects the accuracy of the generalized finite difference approximation.

Because there are no other nodes outside the boundary, if no virtual nodes are added, the node distribution quality in the boundary-node influence domain will be low, i.e., the center of gravity of the local point cloud of the boundary node will deviate greatly from the location of the considered boundary node. The accuracy of the generalized finite difference approximation in Eq. (13) will be reduced at this point. As a result, in GFDM, some studies [34, 35] have identified that a virtual node can be introduced outside the boundary where the boundary node is located to increase the node-distribution quality, and as a result, improve the accuracy of generalized finite difference expressions.

As shown in Fig. 1(a), for the boundary node (marked in red solid points), the nodes (marked in black solid points) in the influence domain (marked in gray area) of the boundary node are all on one side of the tangent line (denoted by the red dotted line) at this boundary node, so the GFDM approximation accuracy of spatial derivatives at this boundary node by using these nodes is low, resulting in low overall calculation accuracy. Therefore, the calculation accuracy of the meshless method is sensitive to derivative boundary conditions (such as Neuman boundary condition and Rudin boundary condition). As shown in Fig. 1(b), the accuracy of generalized finite difference operators at the non-corner boundary node meeting derivative boundary conditions is improved by adding the blue virtual node on the other side of the tangent to make the center of gravity of the local point cloud near to the considered boundary node. That is to say, the index set of the boundary node will be extended by the added virtual nodes.

For the corner boundary node with derivative boundary conditions where there is no tangent, the adding of the virtual nodes needs to be divided into two cases:

Case 1: this case is that the corner node is the intersection of two boundaries, that is, the corner node has two boundary conditions. In this case, the two boundary conditions are generally: (I) a first-type boundary condition and a derivative boundary condition; (II) two derivative boundary conditions; For the corner node at this case, the number of added virtual nodes will be equal to the number of derivative boundary conditions of the corner node, and the direction of adding virtual nodes is the normal direction of the boundary on the side where the corner has derivative boundary conditions. As shown in Fig. 1(c), the left boundary of the red corner node is a derivative boundary condition, and the right boundary is the first-type boundary condition, so you only need to add a blue virtual point along the normal direction of the left boundary. This type of boundary node exists in the numerical examples in Section 3, and the method of adding virtual nodes shown in the point cloud discrete diagrams of the corresponding calculation domains is the processing method described here. As shown in Fig. 1(d), the left and right boundaries of the red corner node are two different derivative boundary conditions, so two blue virtual nodes shall be added along the normal direction of the left and right boundaries respectively. The discrete schemes of the two derivative boundary conditions at the boundary node are the equations of the two virtual nodes respectively, ensuring that the global equations are closed.

Case 2: this case is that the boundaries on both sides of the corner node are the same boundary and share the same derivative boundary condition. In this case, the virtual node is added to the angular bisector of the included angle at the corner. As shown in Fig. 1(e), the red corner node is on a boundary, then the blue virtual node is added to the angular bisector represented by the dotted line. However, this case should be rare, because derivative boundary conditions generally contain normal derivatives, but there is no normal vector at the corner node.

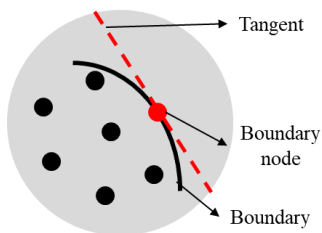

(a) the black nodes in the influence domain of the red non-corner boundary node

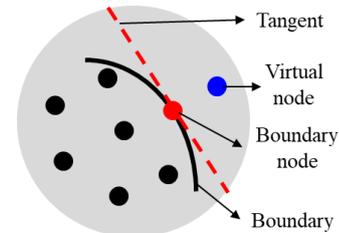

(b) a blue virtual node is added to the influence domain of the red non-corner boundary node

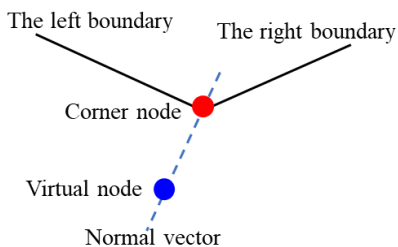

(c) The left boundary of the red corner is the derivative boundary condition, and the right boundary is the first boundary condition in case 1

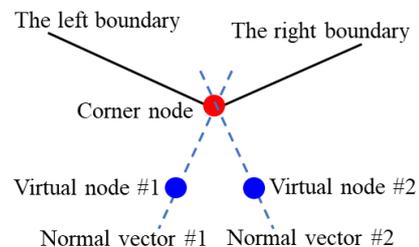

(d) The left boundary and the right boundary of the red corner are derivative boundary conditions in case 1

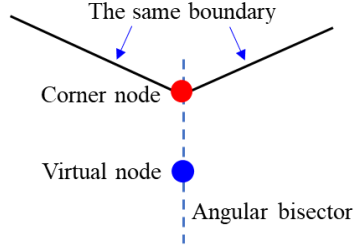

(e) Both sides of the corner node share the same derivative boundary condition in case 2

Fig. 1 Sketch of adding virtual nodes for boundary nodes with derivative boundary conditions

Next, this section takes pressure calculation as an example to illustrate the specific processing details of the derivative boundary conditions.

Assume that there are $n_1$ internal nodes in the computing domain, $n_2$ nodes meeting Dirichlet boundary conditions, and $n_3$ nodes meeting derivative boundary conditions. For the convenience of depiction, assume here that the expressions of the Dirichlet boundary condition and the derivative boundary condition are:

$$p|_A = p_1, \quad \left(hp + l\frac{\partial p}{\partial x}\right)\bigg|_B = q \tag{28}$$

where $A$ and $B$ are the nodes that meet the boundary conditions of the first type and the derivative type respectively, $h$, $l$, and $q$ are coefficients.

For node $A$ that meets the Dirichlet boundary condition, if the sequence number of the node $A$ in all nodes is $a$, the equation corresponding to the boundary condition is:

$$p_a^{n+1} = p_1 \tag{29}$$

For node $B$ that meets derivative boundary conditions, set the serial number of node $B$ in all nodes is $b$, and add the virtual node corresponding to node $B$, denoted as node $C$, because each derivative-boundary-condition node needs to add a corresponding virtual node (Of course, if it is a virtual point shown in Fig. 1(d), two virtual nodes need to be added accordingly. For the convenience of discussion, it is assumed that there is no such virtual node here), the number of nodes in the entire computing domain is $n_1+n_2+n_3+n_3$, if the serial number of the virtual node $C$ in all nodes is $c$, the equation at node $B$ is no longer the equation corresponding to the boundary condition, but the Eq. (30) which is the same as the Eq. (25) of an inner node, and the equation corresponding to the boundary condition is used as the equation corresponding to virtual node $C$. They are:

The discrete equation at node $B$ is

$$\alpha \sum_{j\in \Lambda_b} \left[\frac{k_{bj}}{\mu_{bj}}\left(m_{3j}^b + m_{4j}^b\right)\left(p_j^{n+1} - p_b^{n+1}\right)\right] + q_b^{n+1} = \left[1 - \frac{1-\phi_0}{\phi_0}C_{Temp}\left(T_b^n - T_0\right)\right]C_t \frac{p_b^{n+1} - p_b^n}{\Delta t} \tag{30}$$

The discrete equation at node $C$ is

$$\left(hp + l\frac{\partial p}{\partial x}\right)\bigg|_B = hp_b + l\sum_{j\in \Lambda_b} m_{1j} p_j = q \tag{31}$$

where $n_b$ is the number of nodes in the influence domain of node $B$, $j$ is the serial number of a node in the influence domain of node $B$. Since virtual node $C$ is in the influence domain of node $B$, $c \in \Lambda_b$ holds.

Finally, the linear equations composed of $n_1+n_2+n_3+n_3$ equations can be obtained in the entire computational domain, including $n_1$ Eq. (25), $n_2$ Eq. (29), $n_3$ Eq. (30), and $n_3$ Eq. (31), thus solving the linear equations in a closed manner to obtain the pressure values of all nodes (including $n_1$ inner nodes, $n_2+n_3$ boundary nodes, and $n_3$ virtual nodes) at $n+1$ time step. Combined with Eq. (27) and boundary conditions about temperature (using the same treatment of boundary conditions in this section) calculate the temperature values of all nodes at the $n+1$ time step. Then continue to solve the pressure and temperature distributions at $n+2$ time step.

## 2.5 Application of the upwind GFDM to meshless solution of the convection-diffusion equation and the dissipation error

Different from the diffusion equation, an asymmetric convection term exists in the convection-diffusion equation. When the convection effect is relatively strong, theoretically, it is necessary to use the upwind treatment to discretize the convection term to eliminate the oscillation of the solution caused by the asymmetry of the convection term. As described in Section Introduction, at present, the partial influence domain or the upwind influence domain is mostly used in meshless methods to handle convection-dominated problems. However, due to the possible complexity and changes of the underground flow field, it may be difficult to obtain a stable upwind solution and a general-purpose numerical framework by modifying the influence domain. The upwind GFDM developed in this paper attempts to give a new upwind processing method in the meshless framework. This section will first introduce how to apply the upwind GFDM to conduct a meshless solution of the convection-diffusion equation in theory, and then take the one-dimensional constant-coefficient convection-diffusion equation as an example to illustrate that the discretization of the upwind GFDM for the thermal convection term can be degenerated to the discretization of the convection term by FDM with first-order upwind scheme.

The numerical examples in Sections 3.1, 3.2, and 3.3 demonstrate that the upwind GFDM can realize the effective meshless calculation of single-phase heat and mass transfer problems. Therefore, by constructing a heat and mass transfer problem equivalent to the studied convection-diffusion equation and solving it via the upwind GFDM, the stable meshless solution of the convection-diffusion equation can be obtained.

Assuming that the flow velocity field $\mathbf{V} = (v_x, v_y)$ has been calculated, if it is a steady velocity field, the temperature distribution of the flow field is characterized by the following convection-diffusion equation:

$$\frac{\partial}{\partial x}\left(\lambda_x \frac{\partial T}{\partial x}\right) + \frac{\partial}{\partial y}\left(\lambda_y \frac{\partial T}{\partial y}\right) - v_x \frac{\partial T}{\partial x} - v_y \frac{\partial T}{\partial y} = 0 \tag{32}$$

where $\lambda_x$ and $\lambda_y$ are heat conduction coefficients in $x$ and $y$ directions, and $C$ is the heat capacity.

If it is an unsteady velocity field, the temperature profile of the flow field meets Eq. (33), which is the general form of the convection-diffusion equation.

$$\frac{\partial}{\partial x}\left(\lambda_x \frac{\partial T}{\partial x}\right) + \frac{\partial}{\partial y}\left(\lambda_y \frac{\partial T}{\partial y}\right) - v_x \frac{\partial T}{\partial x} - v_y \frac{\partial T}{\partial y} = C \frac{\partial T}{\partial t} \tag{33}$$

Extend the convection term in Eq. (33) to a second-order derivative term with pressure, that is:

$$-v_x \frac{\partial T}{\partial x} - v_y \frac{\partial T}{\partial y} = -\frac{\partial T v_x}{\partial x} - \frac{\partial T v_y}{\partial y} = \nabla \cdot \left(T \frac{k}{\mu} \nabla p\right) \tag{34}$$

Thus, the convection-diffusion equation in Eq. (34) is extended to the heat and mass transfer coupling problem governed by Eq. (35) and Eq. (36).

$$\frac{\partial}{\partial x}\left(\lambda_x \frac{\partial T}{\partial x}\right) + \frac{\partial}{\partial y}\left(\lambda_y \frac{\partial T}{\partial y}\right) + \nabla \cdot \left(T \frac{k}{\mu} \nabla p\right) = C \frac{\partial T}{\partial t} \tag{35}$$

$$-\frac{k}{\mu} \nabla p = (v_x, v_y) \tag{36}$$

Therefore, when node $i$ is the considered node, the discrete scheme of the convection term can adopt the same discrete scheme as Eq. (22), that is:

$$-v_x \frac{\partial T}{\partial x} - v_y \frac{\partial T}{\partial y} = \nabla \cdot \left(T \frac{k}{\mu} \nabla p\right) = \sum_{j \in \Lambda_i} \left[ T_{ij} \frac{k_{ij}}{\mu_{ij}} (m_{3j} + m_{4j})(p_j - p_i) \right] \tag{37}$$

The nodal temperature values can be calculated by using the sequential coupling or fully-implicit solution of Eqs. (35) and (36), that is, the convection-diffusion equation in Eq. (34) is solved by the upwind GFDM, and the good accuracy and convergence will be illustrated in the numerical example in Section 3.1.

Next, assuming we are studying a 1D problem with a constant flow velocity in the $x$-direction (The flow direction is positive along the $x$-axis, i.e., $v_x > 0$, $v_y = 0$) and homogeneous physical parameters, like the numerical example in Section 3.1. A Cartesian point cloud is used to discretize the rectangular domain. Fig. 2 (a) shows the local point cloud of node 0 in the Cartesian point cloud, i.e. $\Lambda_0 = \{1, 2, 3, 4, 5, 6, 7, 8\}$. May as well suppose $\Delta x = \Delta y$, $r_m = 1.6\Delta x$, then

$$-v_x \frac{\partial T}{\partial x} - v_y \frac{\partial T}{\partial y} = -v_x \frac{\partial T}{\partial x} \tag{38}$$

$$T_{01}=T_{05}=T_{07}=T_1=T_5=T_7, \quad T_{02}=T_{06}=T_{08}=T_2=T_6=T_8, \quad T_0=T_3=T_4$$

where $v_x$ is a constant.

$$p_1 = p_5 = p_7, \quad p_2 = p_6 = p_8, \quad p_0 = p_3 = p_4, \tag{39}$$
$$-(p_1-p_0)=-(p_5-p_0)=-(p_7-p_0)=(p_2-p_0)=(p_6-p_0)=(p_8-p_0)$$

And

$$v_x = -\frac{k}{\mu}\frac{p_1-p_0}{\Delta x} = -\frac{k}{\mu}\frac{p_5-p_0}{\Delta x} = -\frac{k}{\mu}\frac{p_7-p_0}{\Delta x} = \frac{k}{\mu}\frac{p_2-p_0}{\Delta x} = \frac{k}{\mu}\frac{p_6-p_0}{\Delta x} = \frac{k}{\mu}\frac{p_8-p_0}{\Delta x} \tag{40}$$

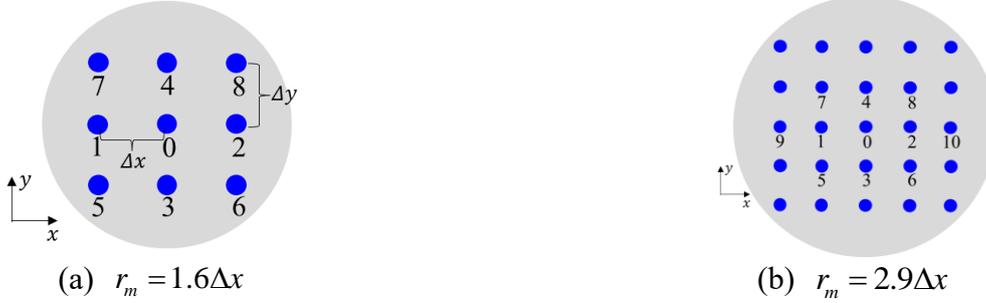

(a) $r_m = 1.6\Delta x$      (b) $r_m = 2.9\Delta x$

Fig. 2 a brief sketch of local point clouds of node 0 with different radius of node influence domain

By introducing Eq. (38), Eq. (39), and Eq. (40) into Eq. (37), we can obtain:

$$-v_x\frac{\partial T}{\partial x}-v_y\frac{\partial T}{\partial y} = \nabla\cdot\left(T\frac{k}{\mu}\nabla p\right) = \sum_{j\in\Lambda_i}\left[T_{ij}\frac{k_{ij}}{\mu_{ij}}(m_{3j}+m_{4j})(p_j-p_i)\right]$$

$$= \frac{k}{\mu}T_1\left[(m_{31}+m_{41})+(m_{35}+m_{45})+(m_{37}+m_{47})\right](p_1-p_0) \tag{41}$$

$$-\frac{k}{\mu}T_2\left[(m_{32}+m_{42})+(m_{36}+m_{46})+(m_{38}+m_{48})\right](p_0-p_2)$$

According to Eq. (12), it can be calculated that:

$$(m_{31},m_{32},m_{33},m_{34},m_{35},m_{36},m_{37},m_{38}) =$$
$$\left(\frac{9.6308\times10^{-1}}{\Delta x^2}, \frac{9.6308\times10^{-1}}{\Delta x^2}, -\frac{3.6917\times10^{-2}}{\Delta x^2}, -\frac{3.6917\times10^{-2}}{\Delta x^2},\right.\tag{42}$$
$$\left.\frac{1.8459\times10^{-2}}{\Delta x^2}, \frac{1.8459\times10^{-2}}{\Delta x^2}, \frac{1.8459\times10^{-2}}{\Delta x^2}, \frac{1.8459\times10^{-2}}{\Delta x^2}\right)$$

Then,

$$\left[(m_{31}+m_{41})+(m_{35}+m_{45})+(m_{37}+m_{47})\right] = \left[(m_{32}+m_{42})+(m_{36}+m_{46})+(m_{38}+m_{48})\right] = \frac{1}{\Delta x^2} \tag{43}$$

Bring Eq. (43) into Eq. (41), it is obtained that:

$$\left(-v_x\frac{\partial T}{\partial x}-v_y\frac{\partial T}{\partial y}\right)\bigg|_0 = \frac{k}{\mu}T_1\frac{(p_1-p_0)}{\Delta x^2} - \frac{k}{\mu}T_0\frac{(p_0-p_2)}{\Delta x^2} = -v_x\frac{T_0-T_1}{\Delta x} \tag{44}$$

As seen in the discrete scheme of the thermal convection term in Eq. (44) derived from the upwind GFDM is exactly the first-order upwind scheme of the convection term in the traditional FDM, which suggests that the upwind GFDM developed in this paper can be regarded as an extension of the traditional upwind FDM in the meshless framework.

Next, we can also analyze the dissipation error of the upwind GFDM and its property. When $r_m = 1.6\Delta x$, the point cloud is as described in Fig. 2 (a). As mentioned earlier, the discrete scheme of the convection term is the first-order upwind scheme in Eq. (44). According to Taylor expansion, it is obtained that:

$$T_1 = T_0 - \Delta x\frac{\partial T}{\partial x}\bigg|_0 + \frac{\Delta x^2}{2}\frac{\partial^2 T}{\partial x^2}\bigg|_0 + O(\Delta x^3) \tag{45}$$

Then, it is obtained that

$$-v_x \frac{\partial T}{\partial x}\Big|_0 = -v_x \frac{T_0 - T_1}{\Delta x} - v_x \frac{\Delta x}{2} \frac{\partial^2 T}{\partial x^2}\Big|_0 + O(\Delta x^2) \tag{46}$$

Thus, the dissipation error of Eq. (44) exists and is:

$$|err_1| = \frac{\Delta x}{2}\left|v_x \frac{\partial^2 T}{\partial x^2}\Big|_0\right| \tag{47}$$

As shown in Fig. 2 (b), when $r_m$ increases to $2.9\Delta x$, without specific calculation, it is expected to obtain a generalized difference expression similar to Eq. (44), that is,

$$\left(-v_x \frac{\partial T}{\partial x} - v_y \frac{\partial T}{\partial y}\right)\Big|_0 = -v_x\left[a\frac{T_0 - T_1}{\Delta x} + (1-a)\frac{T_0 - T_9}{2\Delta x}\right] \tag{48}$$

where $0 < a < 1$.

According to Taylor expansion, it is obtained that:

$$T_9 = T_0 - 2\Delta x \frac{\partial T}{\partial x}\Big|_0 + 2\Delta x^2 \frac{\partial^2 T}{\partial x^2}\Big|_0 + O(\Delta x^3) \tag{49}$$

Then it can be obtained that the dissipation error of Eq. (48) becomes:

$$\begin{aligned}|err_2| &= \left|av_x \frac{\Delta x}{2}\frac{\partial^2 T}{\partial x^2}\Big|_0 + (1-a)v_x \Delta x \frac{\partial^2 T}{\partial x^2}\Big|_0\right| = \left|v_x \frac{\Delta x}{2}\frac{\partial^2 T}{\partial x^2}\Big|_0 + (1-a)v_x \frac{\Delta x}{2}\frac{\partial^2 T}{\partial x^2}\Big|_0\right| \\ &= |err_1| + (1-a)\left|v_x \frac{\Delta x}{2}\frac{\partial^2 T}{\partial x^2}\Big|_0\right| > |err_1|\end{aligned} \tag{50}$$

Similarly, as the radius of the influence domain continues to increase, the dissipation error will continue to increase. The numerical examples in Section 3 will also verify this result.

## 3. Numerical examples

In this section, three numerical examples are designed to test the computational performance and roughly conduct the analysis of error sources of the upwind GFDM, which verify that the method can realize effective calculation of the coupling heat and mass transfer problems.

3.1 A case with a rectangular formation and basic error analysis

In this example, a regular rectangular domain ([0m, 300m]×[0m, 100m]) in Fig. 3(a) is selected. The values of relevant physical parameters are shown in Table 1. The compressibility coefficient, thermal expansion coefficient, and viscosity-temperature coefficient are all 0. The upper and lower boundaries are closed, and the left and right boundaries have constant pressure and temperature. Eq. (51) and Eq. (52) show the specific equations of the physical problem and boundary conditions, thus constructing a stable flow field independent of temperature distribution, to analyze the computational performance of the upwind GFDM to the convection-dominated heat transfer.

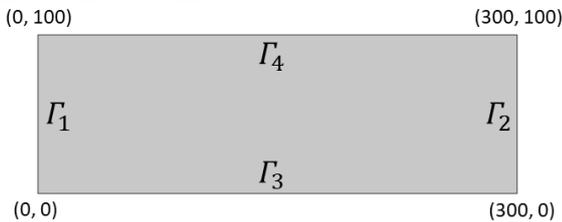

(a) rectangular calculation domain

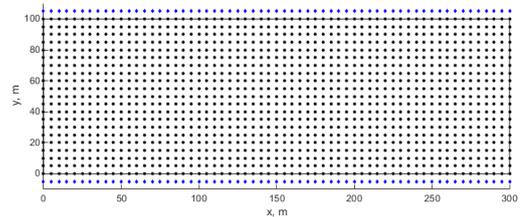

(b) Cartesian point cloud, $\Delta x = 5m$

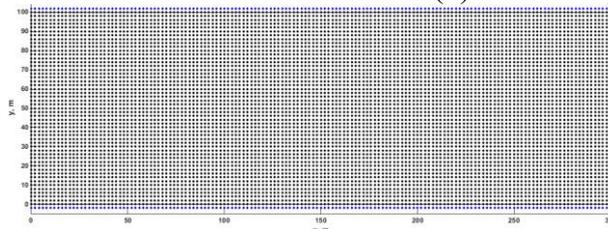

(c) Cartesian point cloud, $\Delta x = 2m$

Fig. 3 Sketches of the calculation domain and the point cloud

Table 1 Values of physical properties of the numerical case

| Parameters | Values |
|---|---|
| Permeability $k$ | 500 mD |
| Heat conduction coefficient of fluid $\lambda_l$ | 0.2 J/s/m/°C |
| Heat conduction coefficient of rock $\lambda_r$ | 3 J/s/m/°C |
| Heat capacity of fluid $C_l$ | $4.2 \times 10^3$ J/KG/°C |
| Heat capacity of rock $C_r$ | 200 J/KG/°C |
| Fluid density $\rho_l$ | 1000 kg/m³ |
| Rock density $\rho_r$ | 2700 kg/m³ |
| Initial porosity $\phi_0$ | 0.3 |
| Compressibility coefficient $C_t$ | 0 MPa⁻¹ |
| Thermal expansion coefficient $C_{Temp}$ | 0 MPa⁻¹ |
| Fluid viscosity at the initial temperature $\mu(T_0)$ | 5 mPa·s |
| Viscosity-temperature coefficient $\alpha_T$ | 0 °C⁻¹ |

$$\nabla \cdot (\nabla p) = 0, \quad p|_{\Gamma_1} = 25, \quad p|_{\Gamma_2} = 10, \quad \frac{\partial p}{\partial y}\bigg|_{\Gamma_3} = 0, \quad \frac{\partial p}{\partial y}\bigg|_{\Gamma_4} = 0 \tag{51}$$

$$186624 \nabla^2 T + 36288000 \nabla \cdot (T \nabla p) = 1638000 \frac{\partial T}{\partial t}, \quad T|_{\Gamma_1} = 40, \quad T|_{\Gamma_2} = 60, \quad \frac{\partial T}{\partial y}\bigg|_{\Gamma_3} = 0, \quad \frac{\partial T}{\partial y}\bigg|_{\Gamma_4} = 0 \tag{52}$$

The analytical solution of Eq. (51) about pressure is: $p = 25 - x/20$.

Thus the pressure gradient is calculated as

$$\nabla p = \left(\frac{\partial p}{\partial x}, \frac{\partial p}{\partial y}\right) = \left(\frac{\partial p}{\partial x}, \frac{\partial p}{\partial y}\right) = \left(-\frac{1}{20}, 0\right) \tag{53}$$

Then Eq. (52) about the heat transfer is rewritten as the following convection-diffusion equation about the temperature.

$$186624 \frac{\partial^2 T}{\partial x^2} - 1814400 \frac{\partial T}{\partial x} = 1638000 \frac{\partial T}{\partial t}, \quad T|_{x=0} = 40, \quad T|_{x=300} = 60 \tag{54}$$

Since Eq. (54) is essentially a 1D heat convection-diffusion problem, this section compares the calculation results of the upwind GFDM, 1D upwind FDM, the reference solution from fine-mesh upwind FDM to test the computational performance of the upwind GFDM and conduct error analysis.

First, the Cartesian point clouds when $\Delta x = \Delta y = 5m$ and $\Delta x = \Delta y = 2m$ shown in Fig. 3(b) and Fig. 3(c) is used to discretize the rectangular domain. Suppose the radius of the node influence domain as $1.6\Delta x$, $2.6\Delta x$, $3.6\Delta x$, $4.6\Delta x$, $5.6\Delta x$, $6.6\Delta x$. Take the time step $\Delta t = 0.5d$, the comparisons of the upwind GFDM results, FDM results, and reference results when $\Delta x = 5m$ and $\Delta x = 2m$ are shown in Fig. 4. Fig. 5 shows some calculated pressure and temperature distributions when $\Delta x = 5m$, $\Delta x = 2m$, and different multiples of the radius of the node influence domain to the node spacing. Fig. 6 shows the L₂ relative error vs. the multiple of the radius of node influence domain to node spacing when $\Delta x = 2m$ and $\Delta x = 5m$. Fig. 7 shows the L₂ relative error of GFDM and traditional FDM vs. the node spacing when $\Delta x = 5m$, $r_m = 1.6\Delta x$, and $r_m = 3.6\Delta x$.

As seen in Fig. 5(h) and Fig. 5(i), the temperature transfer leading edge in the GFDM calculation results gradually becomes curved with the increase of the radius of the influence domain. As seen in Fig. 4, Fig. 5 (h), Fig. 5(i), and Fig. 7, when $r_m = 1.6\Delta x$, the calculation results, calculation accuracy, and convergence of GFDM are almost consistent with those of FDM. As seen in Fig. 4 and Fig. 6, with the increase of the radius of the node influence domain, the range of temperature transfer leading edge in GFDM calculation results becomes larger, and the corresponding calculation error becomes larger. However, even when $r_m$ is large and equal to $3.6\Delta x$, Fig. 7 shows that the upwind GFDM has good convergence. When the node spacing is 1m, the calculation error of the upwind GFDM with $r_m = 3.6\Delta x$ is almost the same as that of the upwind GFDM with $r_m = 1.6\Delta x$. These comparisons show that the upwind GFDM has good computational performance.

We analyze that the error sources of the upwind GFDM calculation results are mainly reflected in the following two points:

Type I error: in the study of two-phase flow in porous media based on the upwind GFDM, Rao et al. [34] pointed out that, the quality of node distribution in the node influence domain (which can also be understood as the deviation degree between the center of gravity of local point cloud in the node influence domain and the node) has a great impact on the approximation accuracy of generalized difference expressions about spatial derivatives, and the larger the radius of the influence domain, the more uneven the node distribution in the influence domain of boundary nodes and nodes close to the boundary, that is, the more the center of gravity of the node in the influence domain of the boundary node deviates from the boundary node, the lower the approximation accuracy of generalized difference expressions of spatial derivatives at boundary nodes and nodes close to the boundary. The result of this type of error reflected in this example is to bend the temperature leading edge that should be a straight line in Fig. 5(h) and Fig. 5(i). For specific details, readers can refer to Rao et al. [34].

Type II error: this type of error is the dissipation error analyzed in Sections 2.5. We know that the dissipation error will widen the leading edge in theory. Taking Fig. 8 (a) of Cartesian collocation as an example, when the radius of the influence domain of node 1 is greater than twice the node spacing, the approximation of heat convection between node 1 and node 3 in the discrete scheme of heat convection term in Eq. (49) is:

$$\rho C_l T_{13} \frac{k_{13}}{\mu_{13}} \left( m_{33}^1 + m_{43}^1 \right) \left( p_1 - p_3 \right) \tag{55}$$

Although node 1 is also upstream of node 3, node 2 is closer upstream from node 3. Heat convection is transmitted through fluid flow. The real fluid flow is from node 1 to node 2, and then from node 2 to node 3. Therefore, since the temperature of node 1 is less than that of node 2, taking node 1 as the upstream temperature in the convection terms of node 1 and node 3 in Eq. (49) accelerates the decline of temperature at node 3. When the most leading edge of the thermal convection reaches node 1, due to the mass transfer between node 1 and node 3, the temperature at node 3 will decrease no matter how small the next time step is, resulting in the wider range of the temperature leading edge than the actual case. As shown in Fig. 4 and Fig. 5, the range of the temperature transfer front increases with the increase of the node influence domain. Under the Cartesian point cloud shown in Fig. 6 and the non-Cartesian point cloud shown in Section 3.2 and Section 3.3, the $L_2$ relative error increases with the increase of the radius of the node influence domain, which verifies the theoretical assertion in Sections 2.5 that the larger the radius of the influence domain, the greater the dissipation error of the upwind GFDM. Imagine a more extreme case. As shown in Fig. 8 (b), if the influence domain of the left-boundary node 4 includes node 5 near the right boundary, the temperature of node 5 will decrease even if the next time step is much small, which is unphysical. The existence of the dissipation error indicates that the radius of the influence domain cannot be set too large when applying the upwind GFDM.

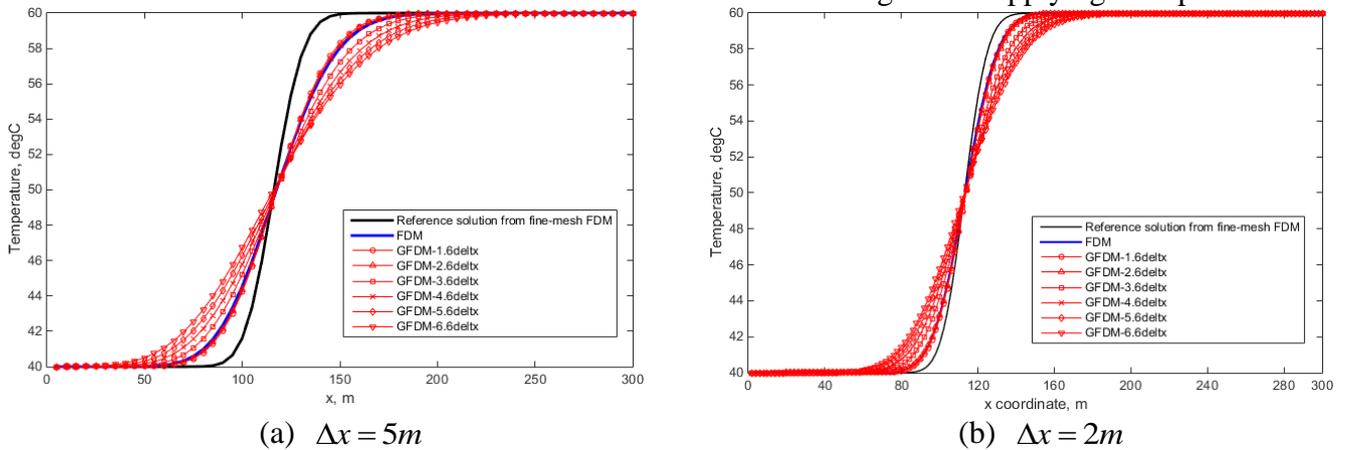

(a) $\Delta x = 5m$  (b) $\Delta x = 2m$

Fig. 4 Comparison of calculation results at section y= 50m from the upwind GFDM and the upwind FDM when $\Delta x = 2m$ and $\Delta x = 5m$

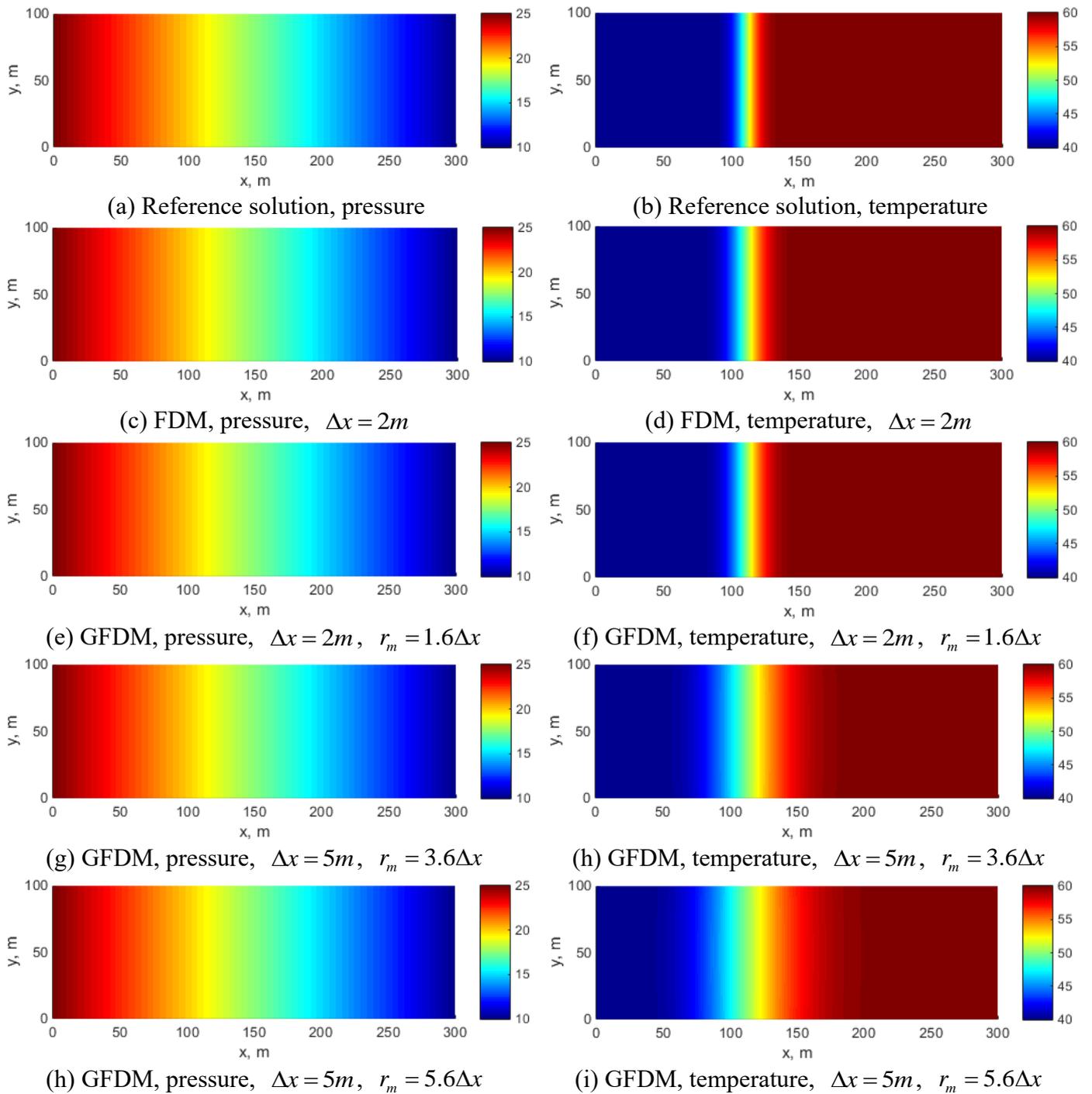

(a) Reference solution, pressure
(b) Reference solution, temperature
(c) FDM, pressure, $\Delta x = 2m$
(d) FDM, temperature, $\Delta x = 2m$
(e) GFDM, pressure, $\Delta x = 2m$, $r_m = 1.6\Delta x$
(f) GFDM, temperature, $\Delta x = 2m$, $r_m = 1.6\Delta x$
(g) GFDM, pressure, $\Delta x = 5m$, $r_m = 3.6\Delta x$
(h) GFDM, temperature, $\Delta x = 5m$, $r_m = 3.6\Delta x$
(h) GFDM, pressure, $\Delta x = 5m$, $r_m = 5.6\Delta x$
(i) GFDM, temperature, $\Delta x = 5m$, $r_m = 5.6\Delta x$

Fig. 5 Comparisons of calculated pressure and temperature profiles by the upwind GFDM and the upwind FDM

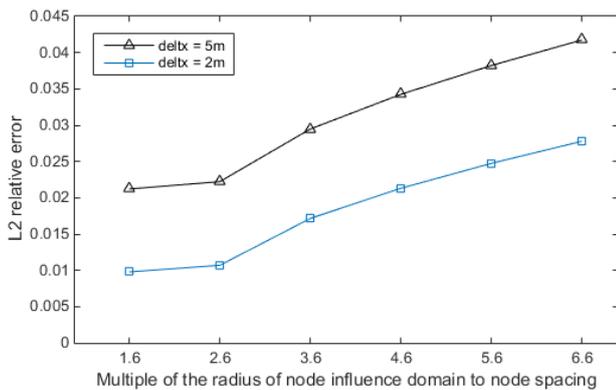

Fig. 6 L$_2$ relative error vs. the multiple of the radius

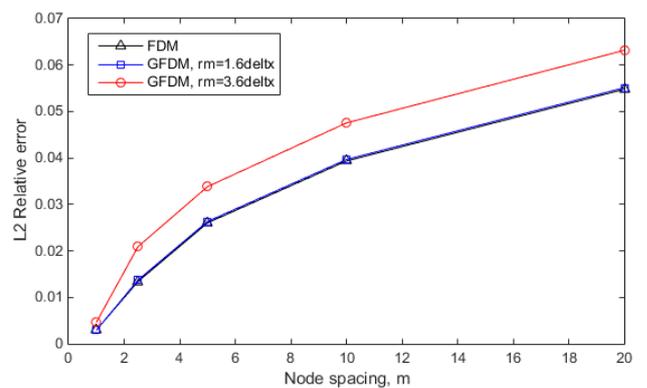

Fig. 7 L$_2$ relative error vs. node spacing when of node influence domain to node spacing when $\Delta x = 2m$ and $\Delta x = 5m$

$\Delta x = 5m$, $r_m = 1.6\Delta x$, and $r_m = 3.6\Delta x$

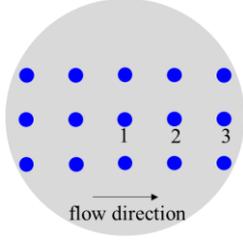
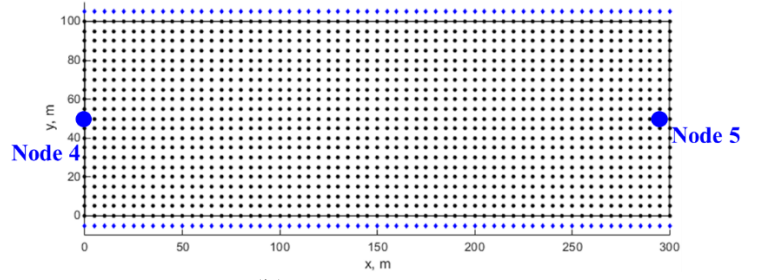

(a) a local point cloud at node 1　　　　　　　　(b) an extreme case

Fig. 8 sketches of some cases for the analysis of Type II error

### 3.2 A case with a polygonal heterogeneous formation

To verify the computational performance of the upwind GFDM when the formation physical parameter is heterogeneous and the domain boundary is irregular, as shown in Fig. 9(a), the example in this section constructs a heterogeneous formation permeability distribution characterized by Eq. (56) and a polygonal calculation domain. Eq. (57) and Eq. (58) show the boundary conditions and initial conditions regarding pressure and temperature respectively. The left boundary and the right boundary are Dirichlet boundary conditions, and the upper and lower boundaries are closed. The values of relevant physical parameters are shown in Table 2.

Fig. 9 shows the fine triangular mesh for FEM, rough triangular mesh for FEM, and the point cloud for GFDM calculation. In this example, the FEM solution based on the fine triangular mesh shown in Fig. 9(b) is used as the reference solution, and the FEM solution based on the rough triangular mesh in Fig. 9(c) and the upwind GFDM solution based on the point cloud in Fig. 9(d) are calculated respectively. In addition, in the calculation process of the upwind GFDM, since boundaries $\Gamma_3$ and $\Gamma_4$ are the second type of boundary conditions, as described in Section 2.4, Fig. 9(d) adds a blue virtual node $6m$ perpendicular to the boundary direction at each blue boundary node to handle the derivative boundary conditions.

$$k = 1200 e^{-x/320} \tag{56}$$

$$p\big|_{\Gamma_1} = 15 MPa,\ p\big|_{\Gamma_2} = 10 MPa,\ p_0 = 10 MPa,\ \frac{\partial p}{\partial y}\bigg|_{\Gamma_3} = 0\, MPa/m,\ \frac{\partial p}{\partial y}\bigg|_{\Gamma_4} = 0\, MPa/m \tag{57}$$

$$T\big|_{\Gamma_1} = 40°C,\ T\big|_{\Gamma_2} = 60°C,\ T_0 = 60°C,\ \frac{\partial T}{\partial y}\bigg|_{\Gamma_3} = 0\,°C/m,\ \frac{\partial T}{\partial y}\bigg|_{\Gamma_4} = 0\,°C/m \tag{58}$$

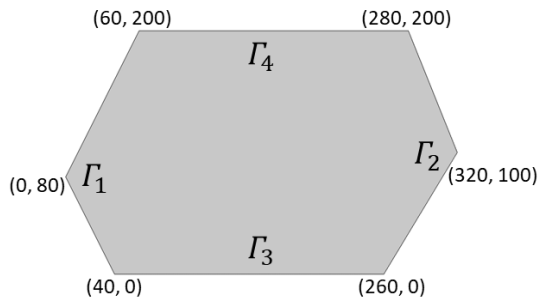
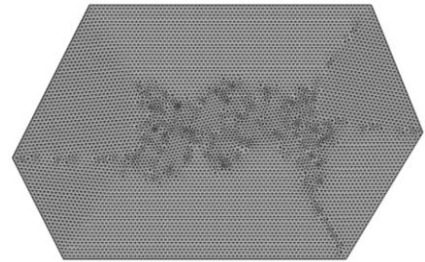

(a) the polygonal calculation domain　　　　　　(b) Fine triangular mesh for reference solution

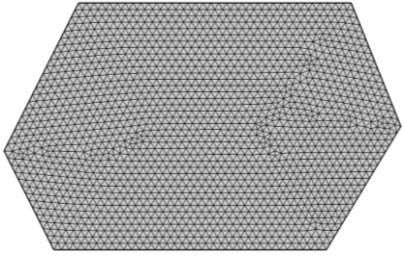
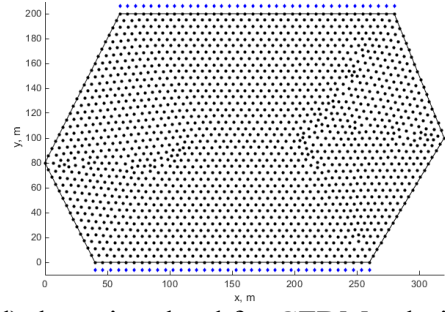

(c) Coarse triangular mesh for FEM solution     (d) the point cloud for GFDM solution

Fig. 9 Sketches of geometric characterizations of the calculation domain

Table 2 Values of physical properties of the numerical case

| Parameters | Values |
| --- | --- |
| Heat conduction coefficient of fluid $\lambda_l$ | 0.2 J/s/m/°C |
| Heat conduction coefficient of rock $\lambda_r$ | 3 J/s/m/°C |
| Heat capacity of fluid $C_l$ | $4.2 \times 10^3$ J/KG/°C |
| Heat capacity of rock $C_r$ | 200 J/KG/°C |
| Fluid density $\rho_l$ | 1000 kg/m$^3$ |
| Rock density $\rho_r$ | 2700 kg/m$^3$ |
| Initial porosity $\phi_0$ | 0.3 |
| Coefficient of compressibility $C_t$ | $1 \times 10^{-5}$ MPa$^{-1}$ |
| Coefficient of thermal expansion $C_{Temp}$ | $1 \times 10^{-5}$ MPa$^{-1}$ |
| Fluid viscosity at the initial temperature $\mu(T_0)$ | 5 mPa·s |
| Viscosity-temperature coefficient $\alpha_T$ | 0.05 °C$^{-1}$ |

Fig. 10 shows the FEM results and the calculated pressure and oil saturation distribution at different influence domain radii (including 15m, 40m, and 65m). Fig. 11 compares the L$_2$ relative errors under different influence domain radii. As seen in Fig. 10 and Fig. 11, (i) for calculated pressure profiles, when $r_m \leq 30m$, the L$_2$ error of the upwind GFDM is lower than that of FEM. (ii) for calculated temperature profiles, the larger the influence domain radius is, the larger the type II error is, and the lower the calculation accuracy of temperature is. But when $r_m = 15m$, the calculation result of the upwind GFDM for temperature distribution is very close to that of FEM, and the relative error of the upwind GFDM for temperature calculation is 0.612%, which is just slightly larger than 0.460% of FEM. It is demonstrated that the upwind GFDM can realize the coupling calculation of mass and heat transfer in the heterogeneous formation, and the calculation error is mainly reflected in the convection-dominated temperature profiles which verifies the error analysis in Section 3.1. It also explains that the treatment of the heterogeneity of physical parameters in Eq. (15) is reasonable, which facilitates the practical application of the upwind GFDM in real complex problems.

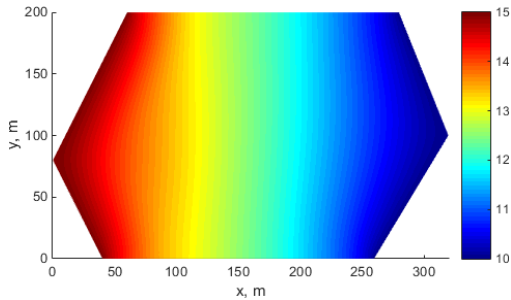
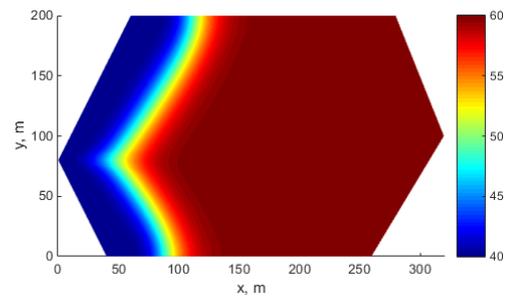

(a) pressure, reference solution     (b) temperature, reference solution

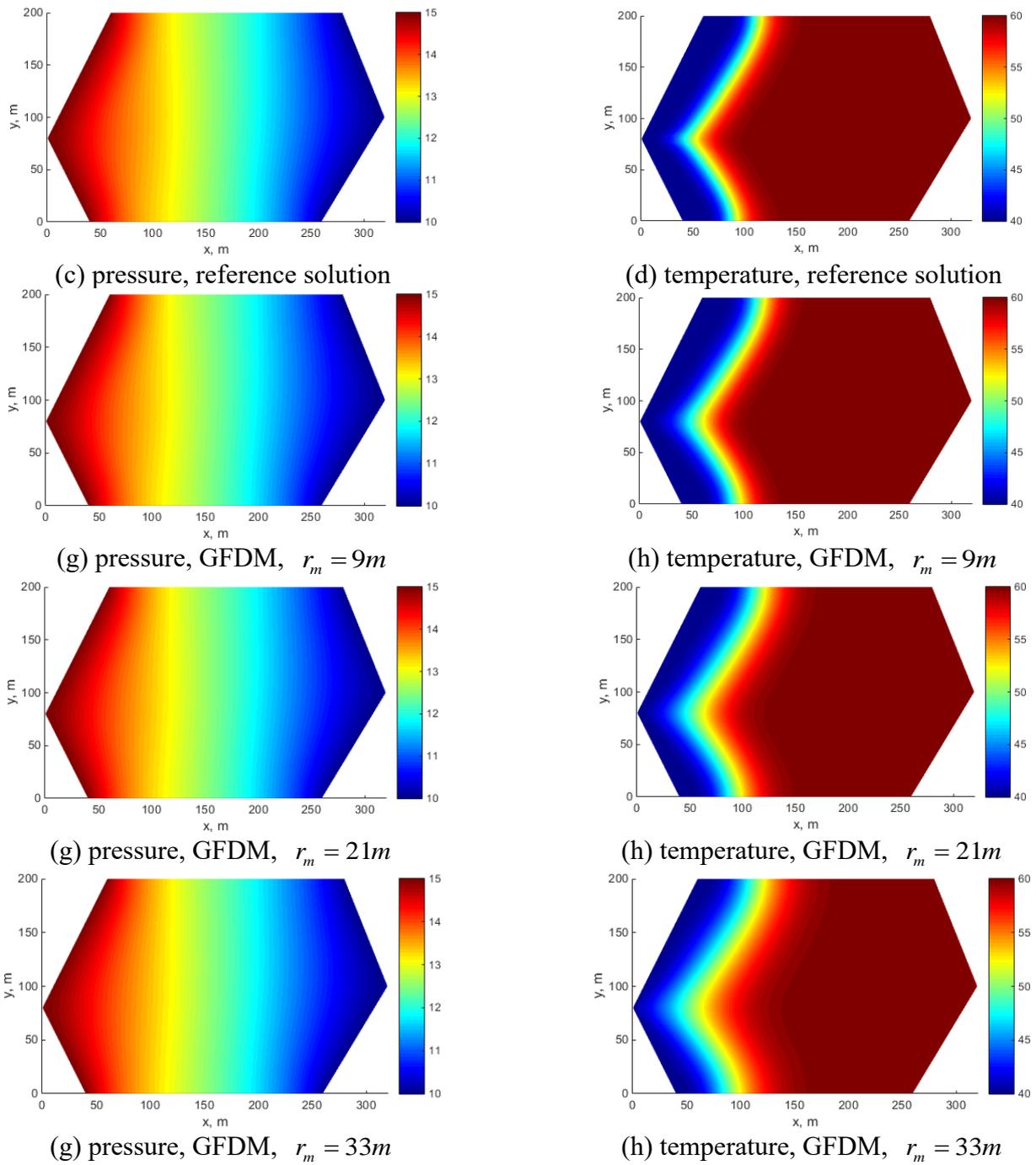

Fig. 10 Comparisons of calculated pressure and temperature profiles by different methods

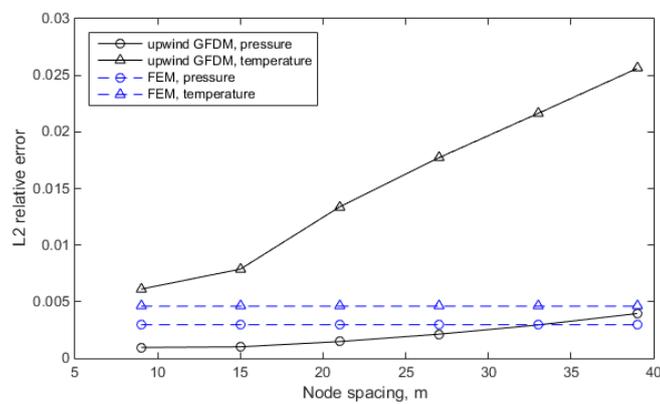

Fig. 11 L$_2$ relative error versus node spacing in example 3

## 3.3 A case with a complex-boundary formation

As shown in Fig. 12(a), the example in this section constructs a more complex calculation domain boundary, in which the blue line segments represent the upper boundary and the lower boundary respectively. Except that the reservoir permeability becomes 1000mD, other physical parameters and boundary conditions are the same as those in Section 3.2. Fig. 12(b) shows the extremely high-density triangular mesh used to obtain the FEM reference solution. Fig. 12(c), (d), (e), and (f) show two different-density triangular meshes for FEM calculation and two different-density point clouds for GFDM calculation respectively. Since the numerical examples and error analysis in Sections 3.1 and 3.2 have demonstrated that the increase of the radius of the influence domain will reduce the calculation accuracy of temperature profiles, when the upwind GFDM calculation is based on the two point clouds in Fig. 12(e) and (f), because the average node spacings in the two point clouds are 12m and 6m respectively, the radius of the node influence domain is taken as 18m and 9m respectively. Fig. 13 and Fig. 14 compare the pressure and temperature distributions at the 40th and 100th days calculated by FEM and upwind GFDM under different-density meshes/point clouds, as well as the reference solutions based on a high-density triangular mesh. Table 3 lists the $L_2$ relative errors of FEM results and upwind GFDM results. It can be seen that when using relatively rough triangulation #1 and point cloud #1, the calculation error of upwind GFDM results for pressure distribution is slightly greater than that of FEM results. The calculation error of the upwind GFDM for temperature distribution at 100 days is nearly twice that of FEM. When using relatively fine triangulation #2 and point cloud #2, the calculation error of upwind GFDM for pressure distribution has rapidly decreased to only half of that of FEM, and the calculation error of temperature has become very close to that of FEM. The comparisons of the calculation errors for pressure profiles in this section and Section 3.2 jointly indicate that the upwind GFDM can handle the pressure diffusion equation more effectively than FEM in a certain radius of the node influence domain and realize higher-accuracy pressure calculation. Overall, the results in this section show that the upwind GFDM can realize the effective calculation of mass and heat transfer in the calculation domain with complex geometry, and implies the good convergence analyzed in Section 3.1.

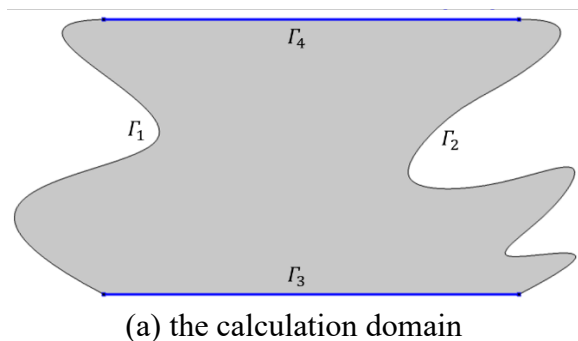

(a) the calculation domain

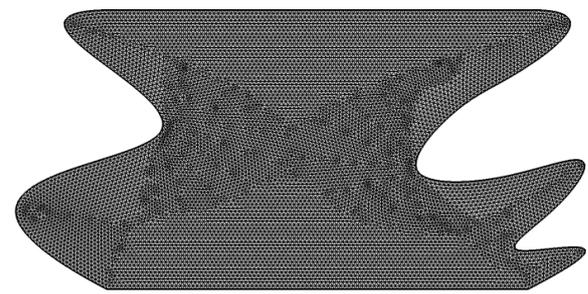

(b) fine triangulation for reference solution

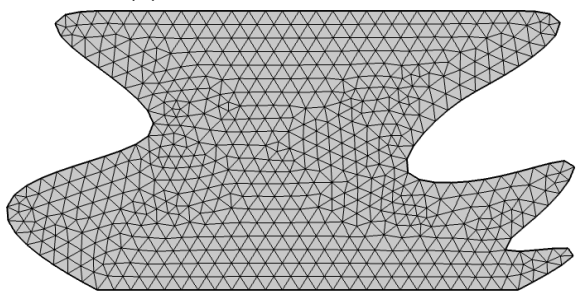

(c) triangulation #1 for FEM

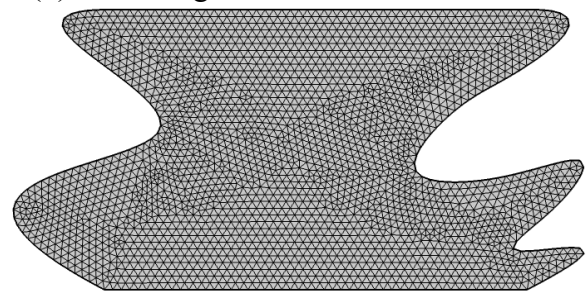

(d) triangulation #2 for FEM

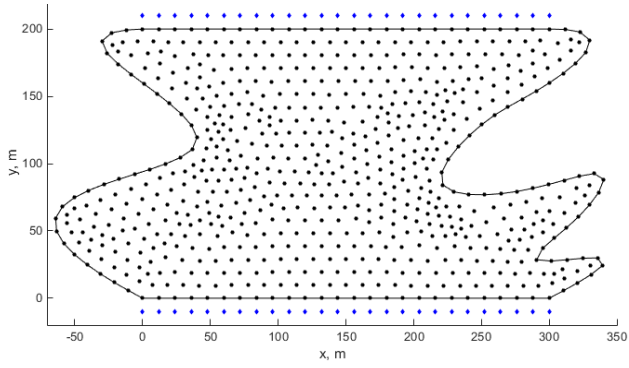
(e) point cloud #1 for FEM

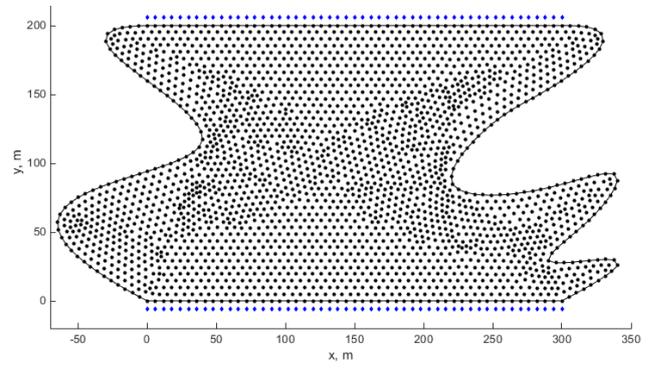
(f) point cloud #2 for FEM

Fig. 12 sketches of the calculation domain, triangulations, and point clouds

Table. 3 $L_2$ relative errors of FEM results and upwind GFDM results

|  | Pressure, 40$^{th}$ day | Temperature, 40$^{th}$ day | Pressure, 100$^{th}$ day | Temperature, 100$^{th}$ day |
| --- | --- | --- | --- | --- |
| FEM (triangulation #1) | 0.183% | 1.805% | 0.174% | 1.487% |
| GFDM (point cloud #1) | 0.214% | 2.434% | 0.231% | 2.871% |
| FEM (triangulation #2) | 0.309% | 1.514% | 0.322% | 1.393% |
| GFDM (point cloud #2) | 0.138% | 1.366% | 0.164% | 1.893% |

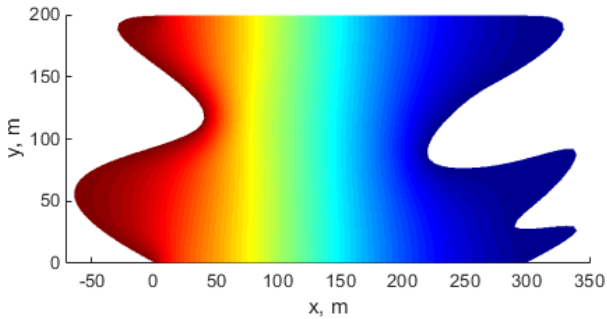
(a) FEM solution, pressure

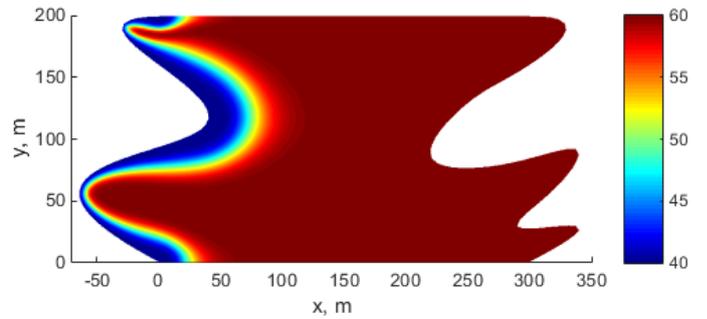
(b) FEM solution, temperature

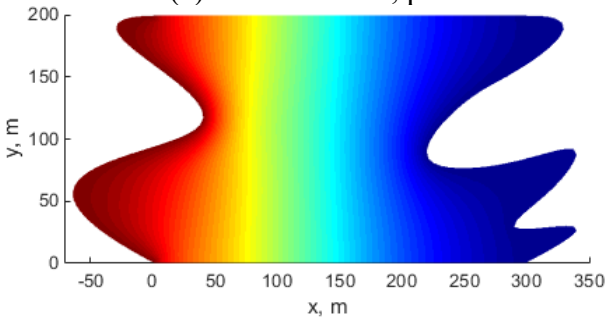
(c) FEM solution with triangulation #2, pressure

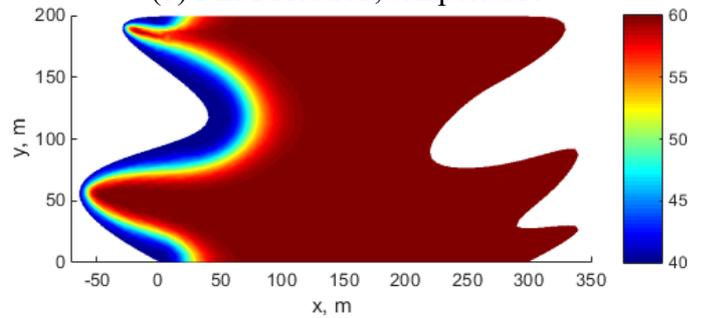
(d) FEM solution with triangulation #2, temperature

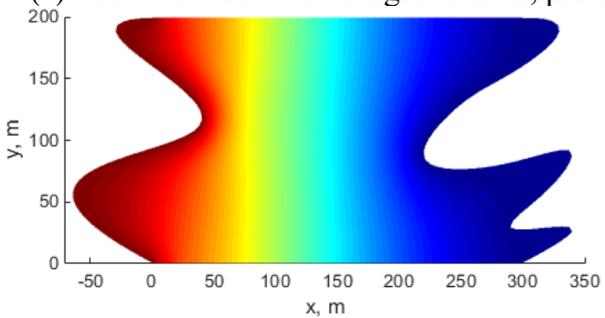
(e) GFDM solution with point cloud #2, pressure

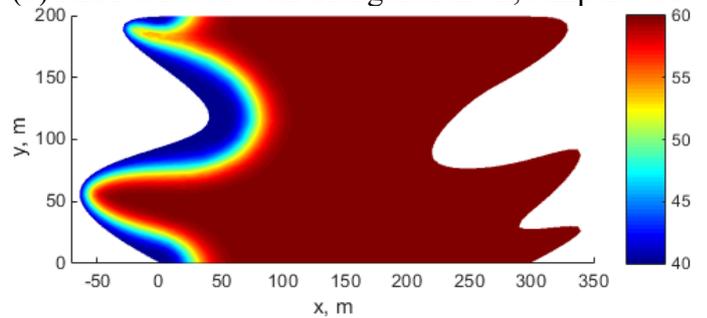
(f) GFDM solution with point cloud #2, temperature

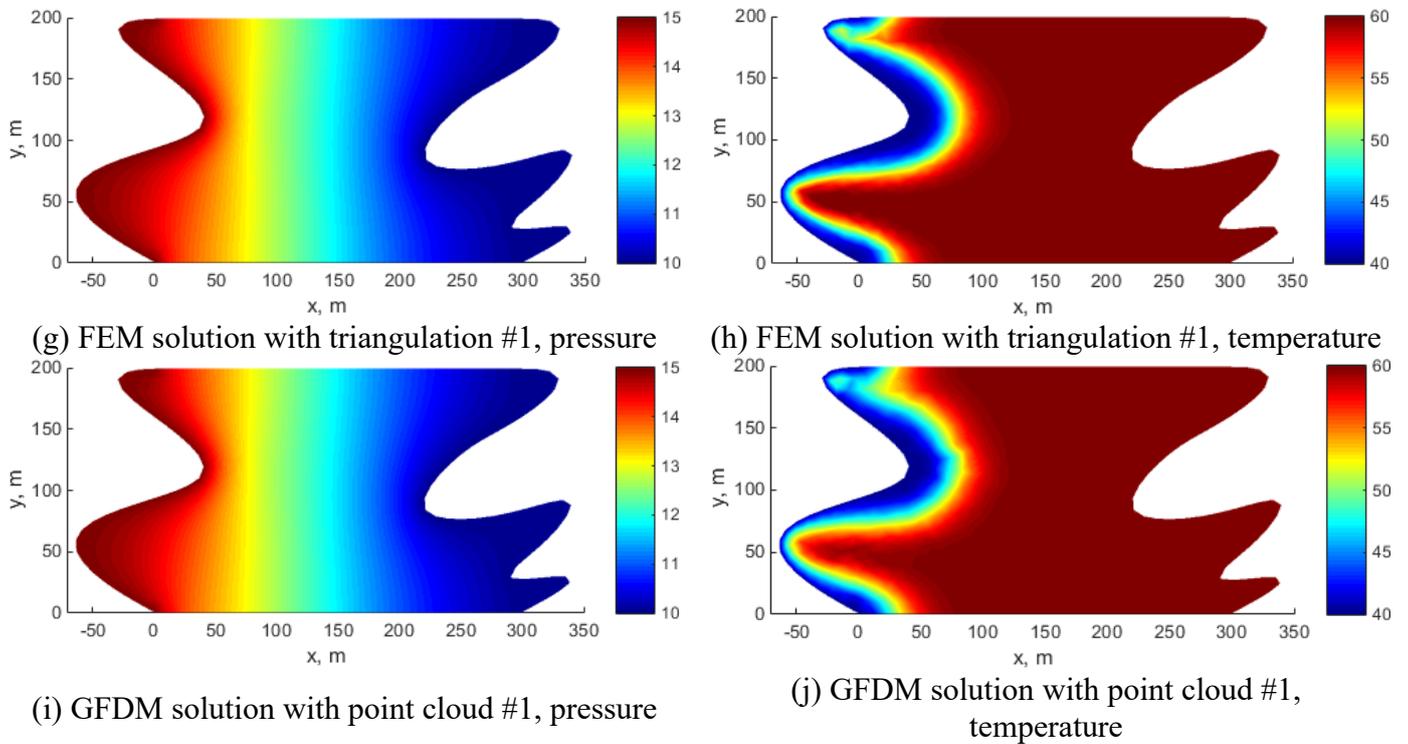

(g) FEM solution with triangulation #1, pressure  (h) FEM solution with triangulation #1, temperature

(i) GFDM solution with point cloud #1, pressure  (j) GFDM solution with point cloud #1, temperature

Fig. 13 Comparisons of calculated pressure and temperature profiles at the 40$^{th}$ day by reference solution, FEM, and the upwind GFDM

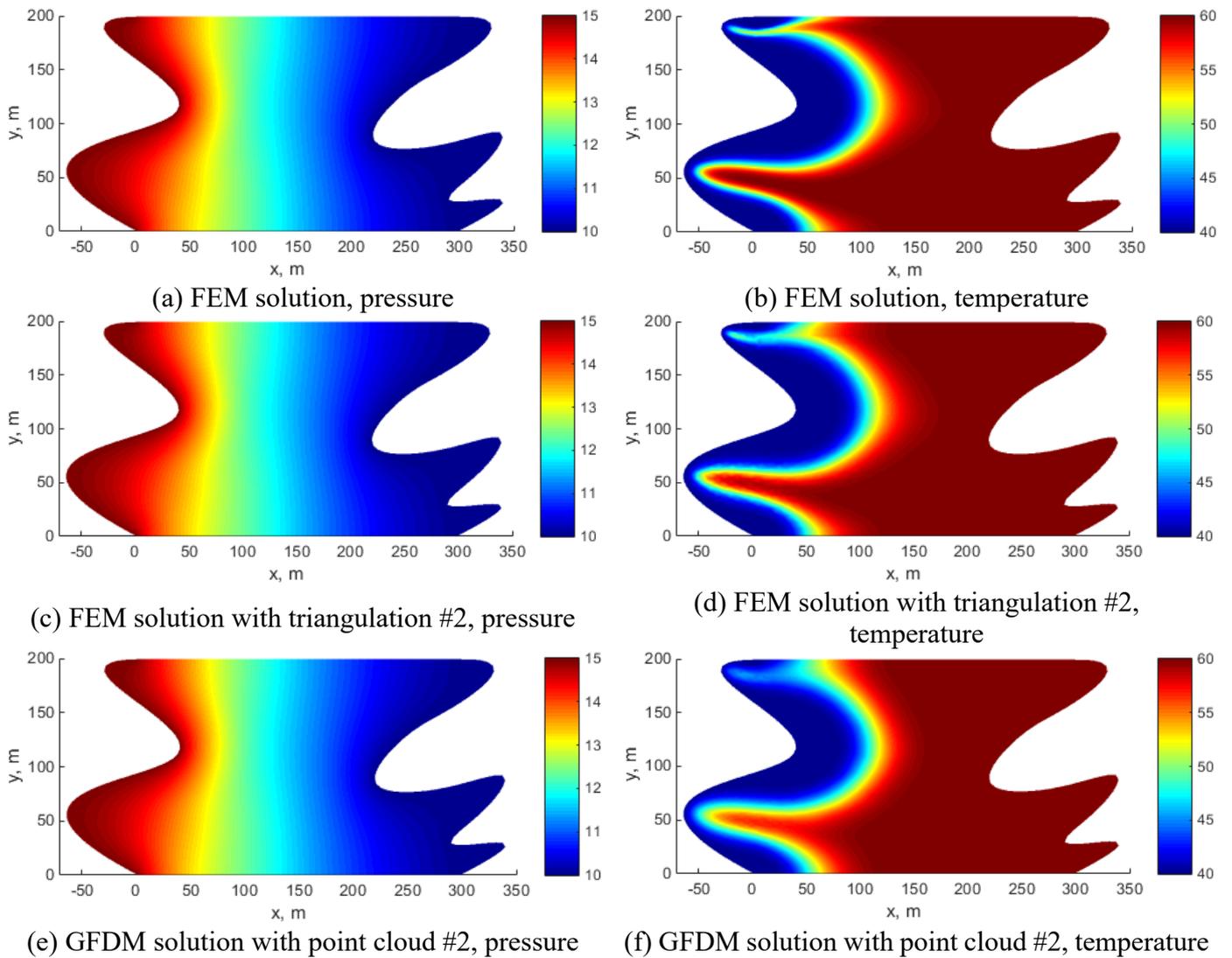

(a) FEM solution, pressure  (b) FEM solution, temperature

(c) FEM solution with triangulation #2, pressure  (d) FEM solution with triangulation #2, temperature

(e) GFDM solution with point cloud #2, pressure  (f) GFDM solution with point cloud #2, temperature

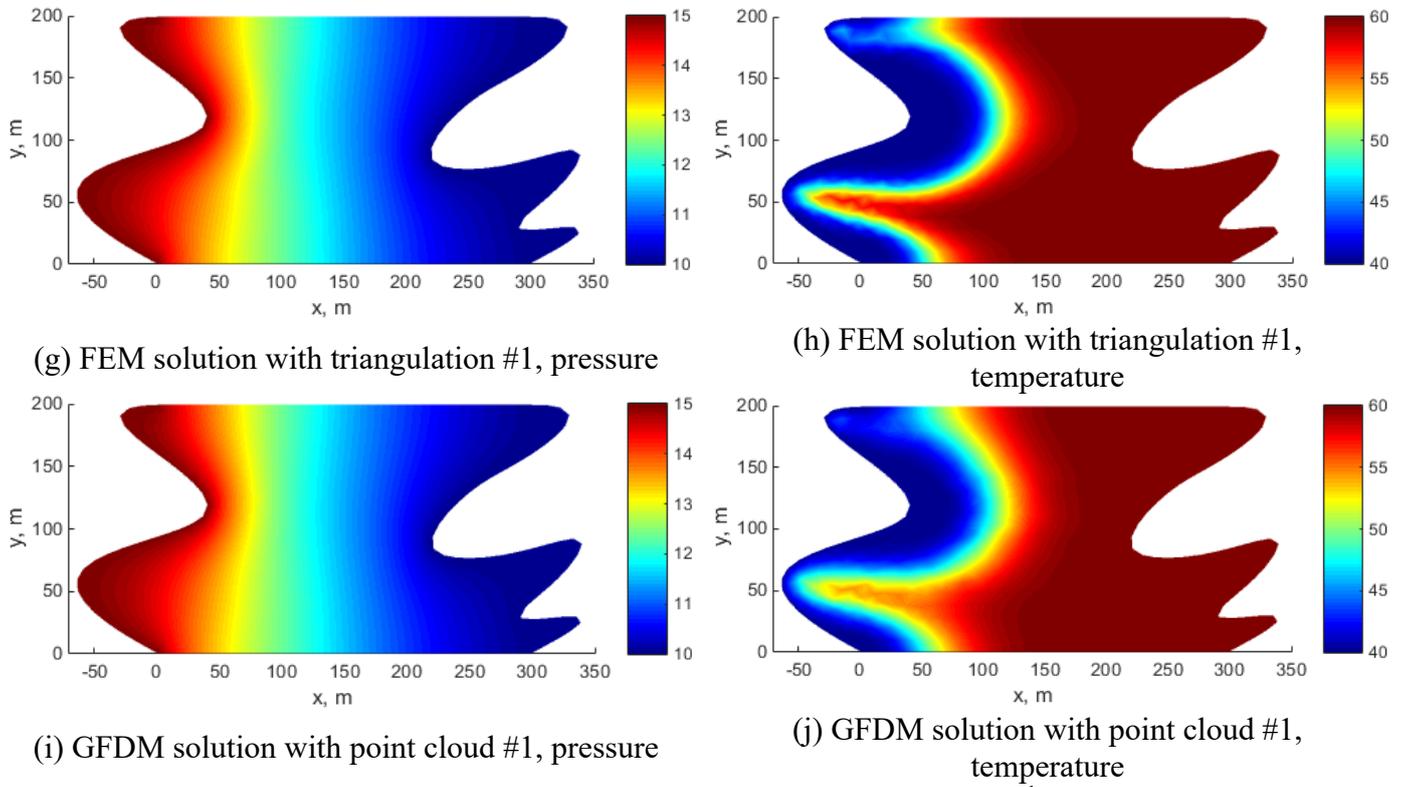

(g) FEM solution with triangulation #1, pressure  (h) FEM solution with triangulation #1, temperature

(i) GFDM solution with point cloud #1, pressure  (j) GFDM solution with point cloud #1, temperature

Fig. 14 Comparisons of calculated pressure and temperature profiles at the 100$^{th}$ day by reference solution, FEM, and the upwind GFDM

## 4. Conclusions and future work

This paper presents an upwind GFDM for heat and mass transfer coupled problems in porous formation. Throughout the whole paper, several key conclusions can be obtained as follows:

(i) This method discretizes the calculation domain by a point cloud rather than the mesh division, due to the topological constraints on the mesh generation being much greater than those on the generation of the point cloud, compared with mesh-based methods, this method can discretize the calculation domain with complex geometry more easily.

(ii) The single point upstream (SPU) scheme commonly used in FDM/FVM-based reservoir simulators is directly introduced to GFDM to form an upwind GFDM. This method is shown able to conduct a meshless solution of the convection-diffusion equation. In addition, take the 1D constant-coefficient convection-diffusion equation as an example, it is proved that the discretization of the upwind GFDM for the heat convection term can be degenerated to the discretization of the convection term by FDM with the first-order upwind scheme, which might indicate the upwind GFDM can be regarded as a meshless extension of the first-order upwind FDM.

(iii) Numerical examples illustrate that the upwind GFDM can realize effective meshless calculation for heat and mass transfer problems in porous media, and obtain a solution of the convection-diffusion equation with a stable upwind effect.

(iv) The upwind GFDM can achieve similar calculation accuracy of temperature profiles but a higher accuracy of pressure profiles compared with the mesh-based methods when the radius of the node influence domain is small, and has good convergence.

(v) This paper analyzes two error sources of the upwind GFDM, and the results of the error analysis and numerical examples both show that the increase of the radius of the node influence domain will increase the calculation error.

The point cloud generation technique for complex 3D computational domains and the computational performance of the upwind GFDM in 3D cases, which we think is possible with significant work in the future. on the other hand, due to the parallel nature of the calculation of the generalized difference operator at each node in the GFDM, the use of parallel computing to form an efficient upwind GFDM-based simulator is an important future work to promote the theoretical research to practice.

## 5. Acknowledgements

Dr. Rao thanks the supports from the National Natural Science Foundation of China (Grant No. 52104017), and the Open Fund of Cooperative Innovation Center of Unconventional Oil and Gas (Ministry of Education & Hubei Province, Yangtze University) (Grant No. UOG2022-14), and Open Fund of Hubei Key Laboratory of Drilling and Production Engineering for Oil and Gas (Yangtze University) (Grant Nos. YQZC202201), and thanks Prof. Hui Zhao and Mr. Wentao Zhan for help in the design of complex-boundary domain.

## 6. Conflict of Interest

The authors declare that they have no conflict of interest.

## 7. Author contributions

Xiang Rao: Conceptualization; Data curation; Methodology; Formal analysis; Funding acquisition; Investigation; Project administration; Software; Validation; Visualization; Writing - original draft; Writing - review & editing.